\documentclass[11pt]{article}

\usepackage{tikz}
\usetikzlibrary{matrix}

\usepackage{amsmath}
\usepackage{amssymb}
\usepackage{latexsym}
\usepackage{enumitem}
\usepackage{mathrsfs}
\usepackage{comment}
\usepackage{color}
\usepackage[normalem]{ulem}

\renewcommand\labelenumi{(\roman{enumi})}
\renewcommand\theenumi\labelenumi

\def\qed{ \ \vrule width.2cm height.2cm depth0cm\smallskip}
\newenvironment{proof}{\noindent {\bf Proof.\/}}{$\qed$\vskip 0.1in}

\newcommand{\ba}{\begin{array}}
\newcommand{\ea}{\end{array}}
\newcommand{\be}{\begin{equation}}
\newcommand{\ee}{\end{equation}}
\newcommand{\bea}{\begin{eqnarray}}
\newcommand{\eea}{\end{eqnarray}}
\newcommand{\beaa}{\begin{eqnarray*}}
\newcommand{\eeaa}{\end{eqnarray*}}

\def\dbE{\mathbb{E}}
\def\dbF{\mathbb{F}}

\def\dbL{\mathbb{L}}

\def\dbN{\mathbb{N}}
\def\dbP{\mathbb{P}}
\def\dbR{\mathbb{R}}
\def\dbS{\mathbb{S}}
\def\dbT{\mathbb{T}}

\def\dbV{\mathbb{V}}

\def\a{\alpha}

\def\d{\delta}
\def\e{\varepsilon}

\def\l{\lambda}

\def\si{\sigma}
\def\t{\tau}
\def\f{\varphi}

\def\o{\omega}

\def\D{\Delta}

\def\O{\Omega}
\def\cA{{\cal A}}

\def\cF{{\cal F}}

\def\cV{{\cal V}}

\def\no{\noindent}

\def\ms{\medskip}
\def\bs{\bigskip}
\def\q{\quad}
\def\qq{\qquad}

\def\pa{\partial}
\def\cd{\cdot}
\def\cds{\cdots}

\def\tr{\hbox{\rm tr}}

\def\qed{ \hfill \vrule width.25cm height.25cm depth0cm\smallskip}

\newcommand{\bas}{\begin{assum}}
\newcommand{\eas}{\end{assum}}

\def\limsup{\mathop{\overline{\rm lim}}}
\def\liminf{\mathop{\underline{\rm lim}}}

\def\pa{\partial}

 \def\cd{\cdot}
\def\cds{\cdots}

\def\tr{\hbox{\rm tr$\,$}}

\def\dis{\displaystyle}

\def\bx{{\bf x}}

\def\1{{\bf 1}}

\def\:{\!:\!}
\def\reff#1{{\rm(\ref{#1})}}
\def \proof{{\noindent \bf Proof\quad}}

at 9pt

\newtheorem{thm}{Theorem}[section]
\newtheorem{lem}[thm]{Lemma}
\newtheorem{cor}[thm]{Corollary}
\newtheorem{prop}[thm]{Proposition}
\newtheorem{rem}[thm]{Remark}
\newtheorem{eg}[thm]{Example}
\newtheorem{defn}[thm]{Definition}
\newtheorem{assum}[thm]{Assumption}

\begin{document}

\title{\bf Dynamic Set Values for Nonzero Sum Games with Multiple Equilibriums}
\author{Zachary {\sc Feinstein}\footnote{Stevens Institute of Technology, School of Business, zfeinste@stevens.edu.}
       \and Birgit  {\sc Rudloff}\footnote{Vienna University of Economics and Business, Institute for Statistics and Mathematics, brudloff@wu.ac.at.}
       \and Jianfeng {\sc Zhang}\footnote{University of Southern California, Department of Mathematics, jianfenz@usc.edu. Research  supported in part by NSF grant DMS-1908665. This author would like to thank Shige Peng for inspiring discussions on the subject. }
}\maketitle

\begin{abstract}
Nonzero sum games typically have multiple Nash equilibriums (or no equilibrium), and unlike the zero sum case, they may have different values at different equilibriums. Instead of focusing on the existence of individual equilibriums, we study the set of values over all equilibriums, which we call the set value of the game. The set value is unique by nature and always exists (with possible value $\emptyset$). Similar to the standard value function in control literature, it enjoys many nice properties such as regularity, stability, and more importantly the dynamic programming principle. There are two main features in order to obtain the dynamic programming principle: (i) we must use closed-loop controls (instead of open-loop controls);  (ii) we must allow for path dependent controls, even if the problem is in a state dependent (Markovian) setting. We shall consider both discrete and  continuous time models with finite time horizon. For the latter we will also provide a duality approach through certain standard PDE (or path dependent PDE), which is quite efficient for numerically computing the set value of the game.  
\end{abstract}
\vspace{5mm}


\vspace{5mm}

\maketitle

\vfill\eject

\section{Introduction}\label{sec:intro}
\setcounter{equation}{0}
{\color{black}
In a standard stochastic control problem, the value function is well-defined and is the unique (viscosity) solution of the associated HJB equation, or the path dependent HJB equation in a path dependent setting. The existence and/or uniqueness of optimal controls often require stronger conditions (typically certain compactness and/or convexity conditions). We remark that the value exists even if there is no optimal control; additionally, when there are multiple optimal controls, they share the same value. Similar results hold for two person zero sum games under the Isaacs condition, where one may study the unique game value without requiring the existence or uniqueness of the equilibriums (saddle points). We refer to the book Mertens, Sorin, \& Zamir \cite{MSZ} for a general exposition of the theory and Possamai, Touzi, \& Zhang \cite[Section 2]{PTZ} for a literature review on continuous time two person zero sum stochastic differential games.  The situation is quite different for nonzero sum stochastic differential games. There have been many works on the existence of Nash equilibriums, by using either PDE method or BSDE method, see e.g.  Bensoussan \& Frehse \cite{BF},  Buckdahn, Cardaliaguet, \& Rainer \cite{BCR}, Cardaliaguet \& Plaskacz \cite{CP}, El-Karoui \& Hamadene \cite{EH}, Friedman \cite{Friedman}, Hamadene \cite{Hamadene},  Hamadene, Lepeltier, \& Peng \cite{HLP}, Hamadene \& Mannucci \cite{HM0}, Hamadene \& Mu \cite{HM1, HM2}, Lin \cite{Lin}, Mannucci \cite{Mannucci1, Mannucci2}, Olsder \cite{Olsder}, Rainer \cite{Rainer}, Sun \& Yong \cite{SY}, Uchida \cite{Uchida}, and Wu \cite{Wu}, to mention a few. We emphasize that, unlike stochastic control problems or zero sum games, in the nonzero sum case different equilibriums could lead to different values,  which makes it difficult to study the game value in a standard manner when there are multiple equilibriums. On the other hand, when there is no equilibrium, it becomes inconvenient even to define  the game value.

We shall define the game value as the set of the values of the game over all equilibriums, which we call the set value of the game. For general set valued analysis, we refer to the book Aubin \& Frankowska \cite{AF}.  With the empty set as a possible set value, both the existence and uniqueness of the set value of the game is always guaranteed by definition. It turns out that this set value behaves benignly as the (real-valued)  value function in stochastic control theory: it enjoys the regularity, stability, and most importantly, the Dynamic Programming Principle (DPP for short) in an appropriate sense. When the set value is a singleton, e.g., in two person zero sum games or in stochastic control problems (a ``game'' with only one player), it reduces to a standard value function (real or vector valued)  and  satisfies a (path dependent) PDE. 
 
Our idea of studying the set value for nonzero sum games follows the line of, among others,  Abreu,  Pearce, \&  Stacchetti \cite{APS} and Sannikov \cite{Sannikov}. {\color{black}The work \cite{APS} considers the set value of an infinitely repeated game in discrete time over all sequential equilibriums. Due to the homogeneousness of the game, its set value is time and state invariant and thus is actually a fixed set, or say a set valued constant. It is shown in \cite{APS} that this set value satisfies the so called factorization and self-generation, which is exactly in the same spirit of our DPP.  The work \cite{Sannikov} considers a similar game, but in continuous time models. The set value is again a fixed set, and the main focuses of \cite{Sannikov} are the characterization and geometric properties of this set as well as their economic implications.} Another highly related work is Cardaliaguet, Quincampoix, \& Saint-Pierre \cite{CQS}, which uses viability theory. The main focus of \cite{CQS} is the numerical approximation for the set of initial states  satisfying some required properties. {\color{black} Our goal is to study standard nonzero sum games in finite time horizon, both in discrete time and in continuous time models, and we shall investigate systematically the dynamic set value of the game  over all Nash equilibriums. }

 In Section \ref{sec:discrete} we study  the discrete time model. Besides establishing the DPP, in the spirit of \cite{APS},  our main contribution is to  show that, even in the state dependent (or say, Markovian) setting, the DPP would fail if one restricts to state dependent equilibriums. Consequently it is necessary to consider path dependent controls in order to have the DPP, which is not the case for stochastic control problems and zero sum games and is due to the non-uniqueness of the values (though the set of values is always unique). While already studied in the literature in various contexts, we also show that DPP would fail if we restrict to Pareto optimal equilibriums and discuss how to choose an ``optimal'' equilibrium by introducing a central planner. Another highly relevant problem, although not discussed in this paper, is to estimate the model parameters with the presence of multiple equilibriums, for which we refer to the survey paper Ho \& Rosen \cite[Section 2]{HR} and the references therein. We shall also remark that, as already observed in Pham \& Zhang \cite{PZ}, through Buckdahn's counterexample for zero sum games, to ensure the DPP for the game value we need to consider closed loop controls rather than open loop controls.

In Section \ref{sec:cont} we study our main object: a continuous time model in a path dependent setting.  
It is in general difficult to study the true equilibriums in this model.  Motivated by \cite{BCR} and \cite[VII.4]{MSZ}, we relax the set value of the game to the limit of the  value sets over all $\e$-equilibriums. Then the set value will be compact and will be nonempty as long as there exist $\e$-equilibriums for all $\e>0$, which is a much weaker requirement than the existence of true equilibriums (see e.g. Frei \& dos Reis \cite{DF} for an example) and is sufficient for practical purpose in most applications. This is exactly in the spirit of the stochastic control problems, where the value is the limit of the values over $\e$-optimal controls. {\color{black}Indeed, for stochastic control problems and zero sum game problems, the (standard) value function corresponds to this relaxed set value, not the original one from true equilibriums when  an optimal control or saddle point does not exist.} We believe this approach of the values could be efficient in more general control/game problems where the optimal control/equilibrium may not exist or is hard to analyze.

Our next result is the regularity (sensitivity with respect to the state process) and stability (sensitivity with respect to the coefficients) of the set value, under mild regularity assumptions on the coefficients. These results have fundamental importance in applications. As a consequence we obtain the measurability of the set value in terms of the state. Our result is in the direction of Feinstein \cite{Feinstein}, except that \cite{Feinstein} studies the set of the equilibriums instead of the values. 

The main result of this paper is the DPP for the set value, which can  be viewed as a type of time consistency and justifies that the set value is an appropriate object for our dynamic model.  While natural in light of its counterpart in the discrete model, the result is much more involved in the continuous time model and requires several approximations. The pathwise setting adds the technical difficulty.  As already observed in Section \ref{sec:discrete}, the pathwise structure is intrinsically needed even in the state dependent setting.

Finally we provide a duality result, motivated by Ma \& Yong \cite{MY} and Karnam, Ma, \& Zhang \cite{KMZ}, which is in the same spirit of the level set approach, see e.g. Barles, Soner, \& Souganidis \cite{BSS}.  We introduce an auxiliary control problem on an enlarged state space, where the additional state  corresponds to the possible values of the game. The value function of the new control problem is a viscosity solution of a standard path dependent HJB equation, for which we refer to Ekren, Touzi, \& Zhang \cite{ETZ1, ETZ2}, and Ren, Touzi, \& Zhang \cite{RTZ}.  Then the set value of the game is characterized as the nodal set of this new value function. This approach is related to the viability approach in \cite{CQS} and is quite efficient in terms of numerical computation of the set value.  
}

\section{The discrete  model}
\label{sec:discrete}
\setcounter{equation}{0}

{\color{black} In this section we study a discrete model with finite time horizon, introduced in Section \ref{sect-setting}. The DPP for the set value is similar to Abreu,  Pearce, \&  Stacchetti  \cite{APS} and is presented in Section \ref{sec:discreteDPP}. The results in Section \ref{sect-discrete-state} concerning the state dependent case is new, to the best of our knowledge. The observations in Sections \ref{sect-static}, \ref{sect-Pareto}, and \ref{sect-optimal} are interesting but not surprising in the game literature. We nevertheless present them here because the same properties hold in the continuous time model in the next section, but it is easier for the readers to include them  in this section.  

\subsection{A static game}
\label{sect-static}
In this subsection we consider a simple static game with $N$ players, and present some basic observations about Nash equilibriums. 
Player $i$'s control takes values in a {\color{black}Borel} measurable set $A_i$ {\color{black}in some arbitrary topological space}. For $a=(a_1,\cds, a_N)\in A := A_1 \times \cds \times A_N$, $J_i(a)$ is the player $i$'s  cost function  she seeks to minimize, and  $J := (J_1,\cds, J_N) : A \to \dbR^N$. We say $a^*\in A$ is a Nash equilibrium if:
\beaa
J_i(a^*) \le J_i(a^{*,-i}, a_i) \q \mbox{for all}~a_i\in A_i,
\eeaa
where  $(a^{*,-i}, a_i)$ is the same as $a^*$ except that its $i$-th component is replaced by $a_i$.

Note that there might be multiple equilibriums or no equilibriums.  We emphasize that  the non-zero sum game  could have different values  $J(a^*)$ at different equilibriums $a^*$, as we see in Example \ref{eg-multiple} below. We thus introduce the set value of the game:
\beaa
\dbV :=\big\{J(a^*): \mbox{for all equilibriums}~ a^*\big\} \subset \dbR^N.
\eeaa

\begin{eg}
\label{eg-multiple}
Set $N=2$, $A_1 = A_2 = \{0, 1\}$, and $J(a)$ as in  Table \ref{tab:multiple} below.  Then  the game has two equilibriums $a^*=(0, 0)$ and $a^* = (1,1)$, and the set value  is $\dbV = \{(0,1), (1,0)\}$. 
\end{eg}

\begin{table}[h]
  \begin{center}
    \begin{tabular}{|l|c|r|} 
     \hline
      $J(a)$ & ${\color{black}a_2 = 0}$ & $a_2 = 1$\\
      \hline
      ${\color{black}a_1=0}$ & $(0,1)$ & $(2,2)$\\
      \hline
      $a_1=1$ & $(3,3)$ & $(1,0)$\\
      \hline
    \end{tabular}
     \caption{\label{tab:multiple} Costs of static non-zero sum game for Example~\ref{eg-multiple}}
  \end{center}
\end{table}

\begin{rem}
\label{rem-discrete-NE}
{\rm The existence of Nash equilibrium is not guaranteed.  However, we emphasize that in this case our set value is still well defined with $\dbV = \emptyset$.  Moreover, our set value is by definition unique, even if there are multiple equilibriums.  
\qed}
\end{rem}

\begin{rem}
\label{rem-discrete-comparison}
{\rm  
\begin{enumerate}
\item\label{rem-discrete-comparison-1} Nash equilibriums may not be Pareto optimal among all controls. Again set $N=2$, $A_1 = A_2 = \{0, 1\}$, and let $J(a)$ be as in the left side of Table \ref{tab:comparison}, then clearly there is a unique equilibrium $a^*=(1,1)$ with value $J(a^*) = (3,3)$. However, we note that $J_i(0,0) =1 < 3 = J_i(a^*)$ for both $i=1,2$. 

\begin{table}[h]
  \begin{center}
    \begin{tabular}{|l|c|r|} 
     \hline
      $J(a)$ & $a_2 = 0$ & $a_2 = 1$\\
      \hline
      $a_1=0$ & $(1,1)$ & $(4, 0)$\\
      \hline
      $a_1=1$ & $(0,4)$ & $(3,3)$\\
      \hline
    \end{tabular}
    \qq\qq \begin{tabular}{|l|c|r|} 
     \hline
      $\tilde J(a)$ & $a_2 = 0$ & $a_2 = 1$\\
      \hline
      $a_1=0$ & $(2,2)$ & $(5, 5)$\\
      \hline
      $a_1=1$ & $(5,5)$ & $(6,6)$\\
      \hline
    \end{tabular}
     \caption{\label{tab:comparison} Costs of static non-zero sum games for Remark~\ref{rem-discrete-comparison}}
  \end{center}
\end{table}

\item\label{rem-discrete-comparison-2} In general  the comparison principle does not hold for the game value. Consider the $\tilde J$ in the right side of  Table \ref{tab:comparison}.  There  is a unique equilibrium $\tilde a^*=(0,0)$ with value $\tilde J(\tilde a^*) = (2,2)$. Note that $J_i(a) < \tilde J_i(a)$ for all $a\in A$ and $i=1,2$, but $J_i(a^*)  =3 > 2=  \tilde J_i(\tilde a^*)$ for both $i=1,2$. 
\qed
\end{enumerate}
}
\end{rem}

\subsection{The set value in a dynamic setting}
\label{sect-setting}
We now consider a dynamic setting. In this section we assume both the time and the state are discrete. Let  $\dbT := \{0,1,\cds, T\}$ denote the set of discrete times, and for each $t\in \dbT$, $\dbS_t$  the set of discrete states at $t$ with $|\dbS_t|<\infty$.  For the reason we will explain in Subsection \ref{sect-discrete-state} below, we shall consider a path dependent setting: $\dbS^\dbT := \big\{\bx = (\bx_0,\cds, \bx_T): \bx_t\in \dbS_t,  t\in \dbT\big\}$.  Set $\O:= \dbS^\dbT$ as the sample space, $\cF:= 2^\O$,  $X_t: \O\to  \dbS_t$  the canonical process: $X_t(\bx) = \bx_t$, and $\dbF = \{\cF_t\}_{0\le t\le T}= \dbF^X$ the natural filtration generated by $X$. Clearly all the functions involved will be $\cF$-measurable.  Throughout this section, all the time dependent functions $\f$ will be required to be adapted in the sense that $\f(t, \bx)$ depends only on $(t, \bx_0,\cds, \bx_t)$. We shall denote
\beaa
\bx =_t \tilde \bx\q\mbox{if}\q \bx_s = \tilde \bx_s~\mbox{for all}~s=0,\cds, t,\q \mbox{and}\q \dbS^\dbT_{t,\bx} := \{\tilde \bx\in \dbS^\dbT: \tilde \bx =_t \bx\}.
\eeaa

There are $N$ players, where the  set of admissible controls $\cA_i$ of the $i$-th player consists of adapted mappings $\a_i:  \dbT\times\dbS^\dbT \to A_i$. Denote  $\cA:= \cA_1\times \cds\times \cA_N$ and $\a := (\a_1, \cds, \a_N)$. For any $(t,\bx, a)\in \dbT\times \dbS^\dbT \times A$, $q(t,\bx, a; \cd): \dbS_{t+1}\to (0,1]$ is a transition probability function: $\sum_{x\in \dbS_{t+1}} q(t,\bx, a; x) = 1$.  Let $\dbP^{t, \bx,\a}$ denote the probability measure such that:
\beaa
\left.\ba{c}
\dis \dbP^{t, \bx, \a}(X =_t \bx) = 1,\q\mbox{and }\\
\dis \dbP^{t, \bx, \a}\big(X_{s+1} = x | X =_s \tilde \bx\big) = q\big(s, \tilde \bx, \a(s, \tilde \bx); x\big)\q\forall~s\ge t,~\tilde \bx \in \dbS^\dbT_{t, \bx},~x\in \dbS_{s+1}.
\ea\right.
\eeaa
Now for $i=1,\cds, N$, let $g_i: \dbS^{\dbT} \to \dbR$ and $f_i :=  \dbT\times \dbS^\dbT\times A_i \to \dbR$ be adapted and measurable in $a_i\in A_i$ (the measurability in $(t, \bx)$ is trivial since the space $\dbT\times \dbS^\dbT$ is finite). The $i$-th player's cost function is defined as:
\beaa
J_i(t, \bx,\a) := \dbE^{\dbP^{t, \bx,\a}} \Big[g_i(X) + \sum_{s=t}^{T-1} f_i(s, X, \a_i(s, X))\Big].
\eeaa
We shall always denote
\beaa
J(t,\bx,\a) := \big(J_1(t,\bx,\a),\cds, J_N(t,\bx,\a)\big) \in \dbR^N.
\eeaa
\begin{defn}
\label{defn-discrete-NE}
Fix $(t,\bx)\in \dbT\times \dbS^\dbT$. We say $\a^*\in \cA$ is a Nash equilibrium of the game at $(t,\bx)$, denoted as $\a^*\in NE(t,\bx)$,  if,  for each $i=1,\cds, N$,
\beaa
J_i(t,\bx, \a^*) \le J_i(t, \bx, \a^{*,-i}, \a_i) \q \mbox{for all}~\a_i\in \cA_i.
\eeaa
\end{defn}

As we saw  in Example \ref{eg-multiple},  the game  could have different values  $J(t,\bx,\a^*)$ at different   equilibriums $\a^*$. Our main object is the following set value over all equilibriums:
\beaa
\dbV(t,\bx) :=\big\{J(t,\bx,\a^*): \a^*\in NE(t,\bx)\big\} \subset \dbR^N,
\eeaa
which is the counterpart of the value function in the standard control literature. As mentioned in Remark \ref{rem-discrete-NE},  $\dbV(t,\bx)$ always exists (with possible value $\emptyset$) and is by nature unique. 

\begin{rem}
\label{rem-discrete-finite}
{\rm {\color{black} For the ease of presentation in this section we restrict to the case $|\dbS_t|<\infty$, but all the results can be easily extended to the case that $\dbS_t$ is countable. When $\dbS_t$ is uncountable, although intuitively the results will still hold true, we will encounter some very subtle measurability issue, as we will see in the next section.  
 }
\qed
}
\end{rem}

\begin{rem}
\label{rem-discrete-zero}
{\rm For two person zero sum games  under the Isaacs condition and other technical conditions, even if there are multiple equilibriums, their values $J$ will always be the same, namely $\dbV(t,\bx) = \{V(t, \bx)\}$ is a singleton, and in the continuous time setting the value function $V$ would satisfy a (path dependent) Isaacs equation. 

{\color{black} We also remark that, by considering mixed strategies, the Isaacs condition will always hold (under very mild conditions), see e.g. Mertens, Sorin, \& Zamir  \cite{MSZ} for discrete time models and Buckdahn, Li,  \& Quincampoix \cite{BLQ} for continuous time models, and hence the set value for these zero sum games is a singleton. It will be interesting to study the set value of nonzero sum games under mixed strategies, which we leave for future research.  }
\qed
}
\end{rem}

We note that, although $\dbS^\dbT$ is finite, unless we assume $A$ is also finite, in general $\dbV(t,\bx)$ may not be finite. The following basic property is interesting in its own right.
\begin{prop}
\label{prop-compact}
If $q$, $f$ are continuous in $a$ and $A$ is compact, then $\dbV(t,\bx)$ is compact.
\end{prop}
\proof  Under our assumption $g(\bx)$ and $f(t, \bx, a)$ are bounded, and thus obviously $\dbV(t,\bx)$ is bounded. Now let $y_n = J(t,\bx, \a_n^*) \in \dbV(t, \bx)$ for some $\a_n^* \in NE(t,\bx)$ and $y_n \to y$. Since  $A$ is compact,  for any $(s, \tilde \bx) \in \dbT\times \dbS^\dbT$,  $\{\a^*_n(s,\tilde \bx)\}_{n\ge 1}$ have a convergent subsequence. Note further that $\dbS^\dbT$ is finite, then without loss of generality we may assume there exists $\a^*\in \cA$ such that $\a^*_n(s,\tilde\bx) \to \a^*(s,\tilde \bx)$ for all $(s, \tilde \bx) \in \dbT\times \dbS^\dbT$. Now for any $i$ and $\a_i\in \cA_i$, we have
\beaa
J_i(t,\bx, \a^*_n) \le J_i(t, \bx, \a_n^{*,-i}, \a_i).
\eeaa
By the continuity of $q$ and $f$ in $\a$, one can easily check that $J_i(t,\bx, \a^*_n) \to J_i(t,\bx, \a^*)$, $ J_i(t, \bx, \a_n^{*,-i}, \a_i) \to J_i(t, \bx, \a^{*,-i}, \a_i)$.  Then $J_i(t,\bx, \a^*) \le J_i(t, \bx, \a^{*,-i}, \a_i)$. This implies  $\a^*\in NE(t,\bx)$, thus $y= J(t,\bx, \a^*)\in \dbV(t,\bx)$.  So $\dbV(t,\bx)$ is closed and hence compact.
\qed

\subsection{Dynamic programming principle for the set value}
\label{sec:discreteDPP}
Given  {\color{black} an $\dbF$-stopping time $\t$  and an $\cF_{\t}$-measurable function $\psi: \dbS^{\dbT} \to \dbR^N$ (namely $\psi(\bx) = \psi(\bx_{\t(\bx)\wedge \cd})$)}, consider the game with terminal time $\t$ and terminal condition $\psi$:
\beaa
J_i(\t, \psi; t, \bx, \a) := \dbE^{\dbP^{t,\bx,\a}}\Big[\psi_i(X) + \sum_{s=t}^{\t -1}  f_i\big(s, X, \a_i(s, X)\big)\Big].
\eeaa
Define the equilibrium at $(\t , \psi; t, \bx)$ in the obvious way and denote its set $NE(\t , \psi; t, \bx)$. Our main result of this section is the following dynamic programming principle.
\begin{thm}
\label{thm-DPP-discrete}
For any $(t, \bx) \in \dbT\times \dbS^\dbT$ and {\color{black}any $\dbF$-stopping time $\t$ with $\t(\bx) > t$}, 
\be
\label{DPP-discrete}
\begin{split}
\dbV(t, \bx) = \Big\{J(&\t , \psi; t, \bx, \a^*):  \mbox{for all $\psi$ and $\a^*$ satisfying} \\
&\psi(\tilde\bx) \in \dbV(\t(\tilde \bx) , \tilde \bx), \forall \tilde\bx\in \dbS^{\dbT}_{t, \bx}, \mbox{and}~ \a^*\in NE(\t , \psi; t, \bx)\Big\}.
\end{split}
\ee
\end{thm}
\proof Let $\tilde \dbV(t, \bx)$ denote the right side of \reff{DPP-discrete}.  

{\it Step 1.} We first prove $\subset$. For any $y = J(t,\bx,\a^*)\in \dbV(t,\bx)$ with $\a^*\in NE(t,\bx)$, denote
\beaa
\psi(\tilde \bx) :=  J(\t(\tilde \bx) , \tilde \bx, \a^*), ~\mbox{for all}~ \tilde \bx \in \dbS^\dbT_{t, \bx}.
\eeaa
Now for any $i$ and $\a_i \in \cA_i$, denote $\tilde \a_i := \a_i\1_{\{s< \t \}} + \a^*_i \1_{\{s\ge \t \}}\in \cA_i$. Then
\beaa
&&J_i(\t , \psi; t, \bx, \a^{*,-i}, \a_i) = \dbE^{\dbP^{t,\bx,\a^{*,-i}, \a_i}}\Big[\psi_i(X) + \sum_{s=t}^{\t -1}  f_i(s, X, \a_i(s, X))\Big]\\
&&\qquad= \dbE^{\dbP^{t,\bx,\a^{*,-i},\tilde \a_i}}\Big[g_i(X) + \sum_{s=t}^{T-1}  f_i(s, X, \tilde \a_i(s, X))\Big] = J_i(t, \bx, \a^{*,-i}, \tilde \a_i).
\eeaa
By setting $\a_i = \a^*_i$ we also have  $J_i(\t , \psi; t, \bx, \a^*)= J_i(t, \bx, \a^*)$. Since $\a^*\in NE(t, \bx)$, then $J_i(\t , \psi; t, \bx, \a^{*,-i}, \a_i)  \ge J_i(\t , \psi; t, \bx, \a^*) $. That is,  $\a^*\in NE(\t , \psi; t, \bx)$. 

Moreover, for any $\tilde \bx \in \dbS^\dbT_{t, \bx}$, denote 
\bea
\label{a1}
\hat \a_i (s, \hat\bx) :=  \a_i  (s, \hat\bx)\1_{\{s\ge \t (\tilde \bx)\}\cap\{\hat\bx ~\!=_{\t(\tilde \bx) } \! ~\tilde \bx\}} + \a^*_i(s, \hat\bx)\1_{(\{s\ge \t (\tilde \bx)\}\cap\{\hat\bx ~\!=_{\t(\tilde \bx) } \! ~\tilde \bx\})^c}\in \cA_i.
\eea
Similarly we have
\beaa
0 &\le& J_i(t, \bx, \a^{*,-i}, \hat \a_i) - J_i(t, \bx, \a^*) \\
&=& \dbP^{t,\bx,\a^*}(X =_{\t(\tilde \bx) } \tilde \bx) \Big[  J_i(\t(\tilde \bx), \tilde \bx,\a^{*,-i}, \a_i)  - \psi_i(\tilde \bx) \Big].
\eeaa
Note that $q>0$ and thus $\dbP^{t,\bx,\a^*}(X =_{\t(\tilde \bx) } \tilde \bx)>0$. This implies that $\a^*\in NE(\t(\tilde \bx) , \tilde \bx)$, then $\psi(\tilde \bx) \in \dbV(\t(\tilde \bx) , \tilde\bx)$. Therefore, it follows from \reff{DPP-discrete} that $y\in \tilde \dbV(t, \bx)$.

{\it Step 2.} On the other hand, let $y=J(\t , \psi; t, \bx, \a^*)\in \tilde \dbV(t, \bx)$ for some desired $\psi$ and $\a^*$.  For each  $\tilde \bx \in \dbS^{\dbT}_{t, \bx}$, we have $\psi(\tilde \bx) \in \dbV(\t(\tilde \bx), \tilde \bx)$ and thus there exists $\a^*_{\tilde \bx}\in NE(\t(\tilde \bx), \tilde \bx)$ such that $\psi(\tilde\bx) = J(\t(\tilde \bx), \tilde \bx, \a^*_{\tilde \bx})$.   Define 
\beaa
\hat \a^*(s, \hat\bx) := \a^*(s, \hat\bx) \1_{\{s < \t(\hat \bx) \}} + \sum_{\tilde \bx\in \dbS^\dbT} \a^*_{\tilde \bx}(s, \hat \bx) \1_{\{s\ge \t(\hat \bx) \}\cap \{\tilde \bx =_{\t(\hat \bx) } \hat \bx\}}\in \cA.
\eeaa
Note that $\t(\tilde \bx) = \t(\hat \bx)$ when $\tilde \bx =_{\t(\hat \bx) } \hat \bx$. Then, for any $i$ and any $\a_i\in \cA_i$, denoting $\tilde \a_i := \a_i  \1_{\{s < \t \}} + \hat\a^*  \1_{\{s \ge \t \}}\in \cA_i$,
\beaa
&&J_i(t, \bx, \hat \a^{*,-i}, \a_i) - J_i(t, \bx, \hat \a^*) \\
&&\qquad= J_i(t, \bx, \hat \a^{*,-i}, \a_i) - J_i(t, \bx, \hat \a^{*,-i}, \tilde \a_i) + J_i(t, \bx, \hat \a^{*,-i}, \tilde \a_i) - J_i(t, \bx, \hat \a^*)\\
&&\qquad= \sum_{\tilde \bx\in \dbS^\dbT} \dbP^{t, \bx, \a^{*,-i}, \a_i} (X =_{\t(\tilde \bx) } \tilde \bx) \Big[J_i(\t(\tilde \bx) , \tilde \bx,  \a_{\tilde \bx}^{*,-i}, \a_i) - J_i(\t(\tilde \bx) , \tilde \bx,  \a_{\tilde \bx}^*)\Big]\\
&&\qquad\qq + J_i(\t , \psi; t, \bx, \a^{*,-i}, \a_i) - J_i(\t , \psi; t, \bx, \a^*)  \\
&&\qquad\ge 0.
\eeaa
This implies  $\hat \a^*\in NE(t,\bx)$, and thus $y= J_i(t, \bx, \hat \a^*)\in \dbV(t,\bx)$.
\qed

\begin{rem}
\label{rem-nondegenerate}
{\rm The condition $q>0$, implying that  $\dbP^{t,\bx, \a}$ are all equivalent for different $\a$, seems crucial in the proof of Theorem \ref{thm-DPP-discrete}.  This condition is also used in \cite{APS} and  is interpreted as that no player can infer the other players' controls through the observed state process. 

 When $q$ is only  nonnegative, we can prove the partial DPP: $\tilde \dbV(t,\bx) \subset \dbV(t, \bx)$,  where $\tilde \dbV(t,\bx)$ again denotes the right side of \reff{DPP-discrete}, {\color{black}and the inclusion could be strict. However, when the measures are singular, it is too strong to require $\psi(\tilde\bx) \in \dbV(\t , \tilde \bx)$ for all $\tilde\bx\in \dbS^{\dbT}_{t, \bx}$. It will be very interesting to see if it is possible to weaken this requirement in an appropriate way so that the DPP will hold true.}
\qed}
\end{rem}

\begin{rem}
\label{rem-open}
{\rm It is crucial that the control is closed loop: $\a = \a(X_\cd)$. If one uses open loop controls, then DPP typically fails even for zero-sum games. See Buckdahn's counterexample  in Pham \& Zhang \cite{PZ} in a continuous time setting, see also Possamai, Touzi, $\&$ Zhang \cite[Remark 4.4(ii)]{PTZ}. {\color{black} Below we present a counterexample in the discrete time setting.}
\qed}
\end{rem}

{\color{black}
We recall that open loop controls do not depend on the state $X$. In this case, the value of $X$, instead of its distribution, will depend on the control.  
\begin{eg}
\label{eg-open} 
Consider a  two player game with open loop controls as follows. Fix a probability space $(\O, \cF, \dbP)$. Set $\dbT := \{t_0, t_1, t_2\}:= \{0,1,2\}$, $\xi_1, \xi_2$ are independent one dimensional random variables with $\dbE[\xi_i] = 0$, $Var(\xi_i)=1$,  the filtration is $\dbF = \{\cF_{t_j}\}_{j=0,1,2}$ with $\cF_{t_0} := \{\emptyset, \O\}$, $\cF_{t_1} := \si(\xi_1)$, $\cF_{t_2} := \si(\xi_1, \xi_2)$, the controls $\a=(\a^1, \a^2)$ are $\dbF$-adapted and take values in $A_1=A_2 := \dbR$, the state process is: for some constant $\si\ge 0$,
\beaa
X^\a_{t_0} :=0,\q X^\a_{t_1} := \a^1_{t_0} + \a^2_{t_0} + \si\xi_1,\q X^\a_{t_2} := [ \a^1_{t_1} + \a^2_{t_1}] X_{t_1} + \si\xi_2,
\eeaa
and the cost functions are: $g_i(x) := -x$, $f_i(t_1, a) := {1\over 2} |a|^2$, $f_i(t_0, a) := 4 |a|^2 + 2a$, that is,
\beaa
J_i(t_0, 0,\a) := \dbE\Big[ {1\over 2} |\a^i_{t_1}|^2 +  4 |\a^i_{t_0}|^2+ 2 \a^i_{t_0} -X^\a_{t_2}\Big],\q i=1,2.
\eeaa
We note that the game is symmetric for the two players. However, DPP fails for this game:
\bea
\label{openDPP}
\left.\ba{c}
\dis\dbV(t_0, 0) = \Big\{\big(-{3\over 2}[|\si|^2+1], ~-{3\over 2}[|\si|^2+1]\big)\Big\},\\
\dis \tilde \dbV(t_0, 0)= \Big\{\big(-[{3\over 2}|\si|^2 + 4], ~-[{3\over 2}|\si|^2 + 4]\Big\}.
\ea\right.
\eea
We note that, when $\si=0$, the above game is deterministic.  

We first show that the two period game has a unique equilibrium: $\a^{*,i}_{t_0} = -{1\over 2}$, $\a^{*,i}_{t_1} = \si \xi_1-1$, $i=1,2$. Then $J_i(t_0, 0, \a^*) = -{3\over 2}[|\si|^2+1]$  and thus we obtain the $\dbV(t_0, 0)$ in \reff{openDPP}. Indeed, assume $\a^*$ is an arbitrary equilibrium. Fix $\a^{*, 2}$. Note that
\beaa
J_1(t_0, 0, \a^1, \a^{*, 2}) = \dbE\Big[ {1\over 2} |\a^1_{t_1}|^2 +  4 |\a^1_{t_0}|^2+ 2 \a^1_{t_0} -[ \a^1_{t_1} +\a^{*, 2}_{t_1}] [\a^1_{t_0} +\a^{*, 2}_{t_0} + \si \xi_1]\Big].
\eeaa
One can easily see that the unique optimal $\a^1_{t_1}$ satisfies: $\a^{*,1}_{t_1} = \a^1_{t_0} + \a^{*, 2}_{t_0} + \si \xi_1$. Then
\beaa
J_1(t_0, 0, \a^1_{t_0}, \a^{*, 1}_{t_1}, \a^{*, 2}) = \dbE\Big[ 4 |\a^1_{t_0}|^2+ 2 \a^1_{t_0}-{1\over 2} [ \a^1_{t_0} +\a^{*, 2}_{t_0} + \si \xi_1]^2 -\a^{*, 2}_{t_1} [\a^1_{t_0} +\a^{*, 2}_{t_0} + \si \xi_1] \Big].
\eeaa
This is strictly convex in $\a^1_{t_0}$. By the  first order condition  we have:
\beaa
0=\dbE\Big[ 8 \a^{*,1}_{t_0}+ 2 -[ \a^{*,1}_{t_0} +\a^{*, 2}_{t_0}+ \si \xi_1] -\a^{*, 2}_{t_1}  \Big]=7\a^{*,1}_{t_0}+2 - \a^{*, 2}_{t_0} - \dbE[\a^{*, 2}_{t_1}].
\eeaa
Similarly, we have $\a^{*,2}_{t_1} = \a^{*,2}_{t_0} + \a^{*, 1}_{t_0} + \si \xi_1$. Then $\dbE[\a^{*,2}_{t_1}] = \a^{*,2}_{t_0} + \a^{*, 1}_{t_0}$, and thus
\beaa
0=7\a^{*,1}_{t_0}+2 - \a^{*, 2}_{t_0} - [\a^{*,2}_{t_0} + \a^{*, 1}_{t_0}] = 6\a^{*,1}_{t_0} - 2\a^{*, 2}_{t_0} +2.
\eeaa
Similarly we have $ 6\a^{*,2}_{t_0} - 2\a^{*, 1}_{t_0} +2=0$. Then one can easily obtain: $ \a^{*,1}_{t_0} = \a^{*, 2}_{t_0} =-{1\over 2}$. This implies that $\a^{*,1}_{t_1} = \a^{*, 2}_{t_1} =  \si \xi_1-1$.

We next compute $\tilde \dbV(t_0, 0)$. Note that
\beaa
J_i(t_1, x, \a_{t_1}) =  \dbE\Big[ {1\over 2} |\a^i_{t_1}|^2-[  \a^1_{t_1} + \a^2_{t_1}] x - \xi_2 \Big]= \dbE\Big[ {1\over 2} |\a^i_{t_1}|^2-[  \a^1_{t_1} + \a^2_{t_1}] x  \Big].
\eeaa
For fixed $x$, one can easily see that the unique equilibrium is $\tilde\a^{*,1}_{t_1} =\tilde\a^{*,2}_{t_1} = x$ (which, for fixed $x$, is deterministic and hence is an open loop control for the game at the second period).  Then $J_i(t_1, x, \tilde\a^*_{t_1}) = -{3\over 2} x^2$ and thus $\dbV(t_1, x) = \{(-{3\over 2} x^2, -{3\over 2} x^2)\}$.  Now consider the game at the first period with terminal $\psi(x) := (-{3\over 2} x^2, -{3\over 2} x^2)$:
\beaa
J_i(t_1, \psi; t_0, 0, \a_{t_0}) = \dbE\Big[4 |\a^i_{t_0}|^2+ 2 \a^i_{t_0} - {3\over 2}[\a^1_{t_0} + \a^2_{t_0} + \si \xi_1]^2  \Big].
\eeaa
By first order conditions we see that the equilibrium satisfies:
\beaa
8 \tilde\a^{*, i}_{t_0} + 2-3 [\tilde\a^{*,1}_{t_0} +\tilde \a^{*,2}_{t_0}]=0,\q i=1,2.
\eeaa
This implies that $\tilde\a^{*,1}_{t_0} = \tilde\a^{*, 2}_{t_0} = -1$,  then $J_i(t_1, \psi; t_0, 0, \tilde\a^*_{t_0}) = -[{3\over 2}|\si|^2 + 4]$.  
 \end{eg}
}

\begin{rem}
\label{rem-symmetric}
{\rm Motivated by the mean field equilibriums, we call an equilibrium $\a^*$ at $(t, \bx)$  symmetric if  $\a^{*,1} = \cds = \a^{*,N}$. Denote
 \beaa
\dbV_{symmetric}(t,\bx) := \{J(t,\bx, \a^*): \mbox{for all symmetric equilibriums $\a^*$}\}.
\eeaa
Then following the same arguments $\dbV_{symmetric}$ also satisfies DPP:  
\beaa
\left.\ba{c}
\dbV_{symmetric}(t, \bx) = \Big\{J(\t, \psi; t, \bx, \a^*): \mbox{for all $\psi$  and $\a^*$ such that $\a^*$ is a symmetric   } \\
\mbox{equilibrium  at $(\t, \psi; t, \bx)$ and $\psi(\tilde\bx) \in \dbV_{symmetric}(\t(\tilde\bx), \tilde \bx)$ for all $\tilde\bx\in \dbS^\dbT_{t, \bx}$}\Big\}.
\ea\right.
\eeaa
\qed}
\end{rem}

\subsection{The state dependent case}
\label{sect-discrete-state}
In this subsection we consider a state dependent (i.e., Markovian) model:  
\bea
\label{State}
q(t,\bx, a; x) = q(t, \bx_t, a; x),\q g (\bx )= g(\bx_T),\q f_i(t, \bx, a ) = f_i(t, \bx_t, a).
\eea
We shall call a function $\f$ on $\dbT\times \dbS^\dbT$ state dependent if $\f(t, \bx) = \f(t, \tilde \bx)$ whenever $\bx_t = \tilde \bx_t$, and in this case it is natural  to abuse the notation and denote it as $\f(t, \bx_t)$.

We first remark that, in this case we may still have path dependent equilibriums, whose value is different from those of state dependent equilibriums.
\begin{eg}
\label{eg-path}
Set $T=3$, $N=2$,   $A_1 = A_2 = \{0, 1\}$,  and $\dbS^\dbT$ takes values as in Figure \ref{fig:path}.
 \begin{figure}[h]
\centering
\begin{tikzpicture}[>=stealth]

 \matrix (tree) [%
matrix of nodes,
minimum size=0.5cm,
column sep=1cm,
row sep=0.2cm,
]
{
                              &  $s_{11}$   &            & $s_{31}$  \\
$\dbS_t$:\q $s_0$&                    & $s_2$ &                 \\
                             &  $s_{10}$    &            & $s_{30}$ \\                  
};
\draw[->] (tree-2-1) -- (tree-1-2);
\draw[->] (tree-2-1) -- (tree-3-2);
\draw[->] (tree-1-2) -- (tree-2-3);
\draw[->] (tree-3-2) -- (tree-2-3);
\draw[->] (tree-2-3) -- (tree-1-4);
\draw[->] (tree-2-3) -- (tree-3-4);
\end{tikzpicture}
\caption{\label{fig:path} States for Example~\ref{eg-path}}
\end{figure}
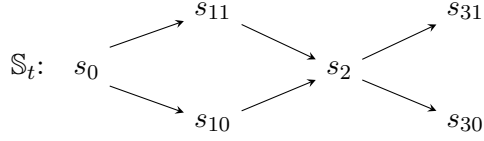

\no That is, $\dbS_0 = \{s_0\}$, $\dbS_1 = \{s_{10}, s_{11}\}$, $\dbS_2 = \{s_2\}$, $\dbS_3=\{s_{30}, s_{31}\}$.  For the first two periods and for $g$, we set
\beaa
f(0,\cd)= f(1,\cd) = 0,\q q(0,\cd) = {1\over 2},\q q(1,\cd) = 1, \q g(s_{30})=(1, 1), \q g(s_{31}) = (0,0),\q
\eeaa
 Then the game at $(0, s_0)$ does not depend on $\a(0,\cd)$ and $\a(1, \cd)$. Indeed,
\bea
\label{pathJ0}
\left.\ba{c}
\dis J(0, s_0, \a) ={1\over 2}\Big[\bar J(\a(2, (s_0, s_{10}, s_2))) + \bar J(\a(2, (s_0, s_{11}, s_2)))\Big], \\
\dis  \mbox{where}\q \bar J_i(a) =  f_i(2, s_2, a_i) + q(2, s_2, a; s_{30}),\q i=1,2.
\ea\right.
\eea
Let us assume that the game for $\bar J(a)$, which corresponds to the last period of the original game, has two equilibriums $a^*$ and $\tilde a^*$. Then we may construct a path dependent equilibrium: noting that $X_0 \equiv s_0$ and $X_2 \equiv s_2$ are deterministic,
\bea
\label{path-equilibrium}
\a^*(2, X) := a^* \1_{\{X_1 = s_{10}\}} + \tilde a^* \1_{\{X_1 = s_{11}\}}.
\eea

 For this purpose, we set  $f(2,s_2, a)$  and $q(2, s_2, a; s_{30})$ for $a\in A$ as in Table \ref{tab:path}. Then by \reff{pathJ0}  we see that $4\bar J$  is the same as Table \ref{tab:multiple}, and thus there are two equilibriums $a^*=(0,0)$ and $\tilde a^*=(1,1)$, with corresponding values $\bar J(a^*) = (0, {1\over 4})$ and $\bar J(\tilde a^*) = ({1\over 4}, 0)$.
 
 \begin{table}[h!]
  \begin{center}
  \begin{tabular}{|l|c|r|} 
     \hline
      $f(2, s_2, a)$ & $a_2 = 0$ & $a_2 = 1$\\
      \hline
      $a_1=0$ & $(-{1\over 4}, 0)$ & $(-{1\over 4}, -{1\over 4})$\\
      \hline
      $a_1=1$ & $(0,0)$ & $(0,-{1\over 4})$\\
      \hline
    \end{tabular}
   \qq\qq
    \begin{tabular}{|l|c|r|} 
     \hline
      $q(2, s_2, a; s_{30})$ & $a_2 = 0$ & $a_2 = 1$\\
      \hline
      $a_1=0$ & ${1\over 4}$ & ${3\over 4}$\\
      \hline
      $a_1=1$ & ${3\over 4}$ & ${1\over 4}$\\
      \hline
    \end{tabular}
      \caption{\label{tab:path} Cost matrices and transition probabilities for Example~\ref{eg-path}}
  \end{center}
\end{table}

We now come back to the original game $J(0, s_0, \a)$. Note that, by \reff{pathJ0}, the only relevant control is $\a(2, (s_0, X_1, s_2))$. If $\a$ is  state dependent, then $\a(2, (s_0, X_1, s_2)) = \a(2, s_2)$ is deterministic. This implies $J(0, s_0, \a) = \bar J(\a(2, s_2))$,  and thus there are only two equilibriums with values $(0, {1\over 4})$ and $({1\over 4}, 0)$. However, we can construct a path dependent equilibrium $\a^*$ by \reff{path-equilibrium}, whose corresponding value is: $J(0, s_0, \a^*)={1\over 2} \bar J(a^*) + {1\over 2} \bar J(\tilde a^*) = ({1\over 8}, {1\over 8})$.
\end{eg}

In view of Example~\ref{eg-path}, nevertheless, $\dbV$ is still state dependent if we restrict to the state dependent model \reff{State}. 
\begin{prop}
\label{prop-state}
Under \reff{State}, $\dbV(t,\bx) = \dbV(t, \bx_t)$ is state dependent.
\end{prop}
\proof Assume $\bx_t = \bx'_t$. For any $\a\in \cA$  and $\tilde \bx' \in \dbS^\dbT_{t, \bx'}$,  introduce $\a'$ by $\a'(s, \tilde \bx') := \a(s, \tilde \bx)$ where $\tilde \bx_s := \bx_s\1_{\{s\le t\}} + \tilde \bx'_s \1_{\{s>t\}}$.  Then one can easily check that $J(t,\bx, \a) = J(t, \bx', \a')$. Such correspondence is one to one and thus it is clear that $\dbV(t,\bx) = \dbV(t,\bx')$.
\qed

From now on, in the state dependent case, we may write the set value as $\dbV(t,x)$. The following DPP is an immediate consequence of Theorem \ref{thm-DPP-discrete}.
\begin{cor}
\label{cor-DPP-state}
Under \reff{State}, for any $(t, \bx) \in \dbT\times \dbS^\dbT$ and  $\dbF$-stopping time $\t$ with $\t(\bx)  \ge t$, 
\beaa
\begin{split}
\dbV(t, x) = \Big\{J(&\t, \psi; t, \bx, \a^*): \mbox{for all $\psi, \a^*, \bx$  such that  $ \bx_{t}=x$, } \\
&\mbox{ $\psi(\tilde\bx) \in \dbV(\t(\tilde \bx), \tilde \bx_{\t(\tilde \bx)})$ for all $\tilde\bx\in \dbS^\dbT_{t, \bx}$,  and $\a^*\in NE(\t, \psi; t, \bx)$ }\Big\}.
\end{split}
\eeaa
\end{cor}

We emphasize that, although our model is state dependent here, the DPP above involves path dependent  $\psi$ and $\a^*$. In fact, if we restrict to state dependent functions $\psi$ and/or $\a^*$, then the DPP may fail, as we explain next. {\color{black}For simplicity, below we consider only deterministic time: $\t \equiv T_0$ for some $T_0 >t$.}

We first investigate the case that $\psi$ is state dependent but $\a^*$ can be still path dependent. In this case by Corollary \ref{cor-DPP-state} the following partial DPP is obvious:
\bea
\label{DPP-state2}
\begin{split}
\dbV(t, x) \supset \Big\{J&(T_0, \psi; t, \bx, \a^*): \mbox{for all  state dependent $\psi$ and $\a^*\in\cA$, $\bx\in \dbS^\dbT$   } \\
&\mbox{such that  }~ \bx_t = x,~ \psi(\tilde x) \in \dbV(T_0, \tilde x), \forall \tilde x\in \dbS_{T_0}, ~ \mbox{and}~ \a^*\in NE(T_0, \psi; t, \bx)\Big\}.
\end{split}
\eea
However, the above inclusion can be strict.
\begin{eg}
\label{eg-psistate}
Consider Example \ref{eg-path} and set $T_0 = 2$.  By Example \ref{eg-path}, we see that
\beaa
\dbV(2, s_2) =\big\{\bar J(a^*), \bar J(\tilde a^*)\big\} = \big\{(0, {1\over 4}), ({1\over 4}, 0)\big\}.
\eeaa
If $\psi$ is state dependent, then there are only two possible functions: $\psi_1(s_2) = (0, {1\over 4})$ and $\psi_2(s_2) = ({1\over 4}, 0)$. Recalling that $f(0,\cd) = f(1,\cd)=0$, then $J(T_0, \psi; 0, s_0, \a) = \psi(s_2)$ for all $\a$.  Thus the right side of \reff{DPP-state2} is $\big\{(0, {1\over 4}), ({1\over 4}, 0)\big\}$. However, by Example  \ref{eg-path} we know that $\dbV(0, s_0)$ contains at least one more value $({1\over 8}, {1\over 8})$.
\end{eg}

We next investigate the case that both $\psi$ and $\a\in \cA$ are state dependent, then obviously $J(t,\bx, \a)$ and $J(T_0, \psi; t, \bx, \a)$ are also state dependent. Define 
\beaa
\left.\ba{c}
\dis\cA_{state} := \{\a\in \cA: ~\mbox{$\a$ is state dependent}\};\\
\dis\dbV_{state}(t, x) := \{J(t,x,\a^*): ~\mbox{$\a^*\in \cA_{state}$ is an equilibrium among all $\a\in \cA_{state}$}\}.
\ea\right.
\eeaa
We emphasize that here all controls are required to be state dependent, in particular, the above $\a^*\in \cA_{state}$ may not be an equilibrium among all controls $\a\in \cA$. Consequently, $\dbV_{state}(t, x)$ may not be a subset of $\dbV(t,x)$. Again,  $\dbV_{state}$ does not satisfy the DPP.

\begin{prop}
\label{prop-DPPstate3}
Under \reff{State}, $\dbV_{state}$ satisfies a partial DPP:
\bea
\label{DPP-state3}
\left.\ba{c}
\dbV_{state}(t, x) \subset \Big\{J(T_0, \psi; t, x, \a^*): \mbox{for all state dependent $\psi$ and $\a^*\in \cA_{state}$  s.t.   } \\
\mbox{$\psi(\tilde x) \in \dbV_{state}(T_0, \tilde x)$, $\forall \tilde x\in \dbS_{T_0}$, and $\a^*$ is an  equilibrium in $\cA_{state}$ at $(T_0, \psi; t, x)$ }\Big\},
\ea\right.
\eea
but the inclusion could be strict. 
\end{prop}
We remark that the inclusions in \reff{DPP-state2} and \reff{DPP-state3} have opposite directions.

\ms
\proof Let $\tilde V_{state}(t, x)$ denote the right side of \reff{DPP-state3}. We shall prove $\dbV_{state} \subset \tilde \dbV_{state}$, and see Example \ref{eg-state} below  that $\dbV_{state} \neq \tilde \dbV_{state}$. 
We follow the arguments in Theorem \ref{thm-DPP-discrete} Step 1
and proceed in two steps.  

{\it Step 1.} Let  $\a^*\in \cA_{state}$ be an equilibrium in $\cA_{state}$ at $(t,x)$.  Denote
\beaa
\psi(\tilde x) :=  J(T_0, \tilde x, \a^*), ~\mbox{for all}~ \tilde x \in \dbS_{T_0}.
\eeaa
For any $i$ and $\a_i \in \cA_{state,i}$, note that $\tilde \a_i := \a_i\1_{\{s< T_0\}} + \a^*_i \1_{\{s\ge T_0\}}$ is also in $\cA_{state, i}$. Then following the same arguments as in Theorem \ref{thm-DPP-discrete} Step 1 
we see that  $\a^*$ is an equilibrium in $\cA_{state}$ at $(T_0, \psi; t, x)$. 

{\it Step 2.} It remains to show that $\psi(\tilde x) \in \dbV_{state}(T_0, \tilde x)$ for all $\tilde x\in \dbS_{T_0}$. That is,
\bea
\label{state-induction}
J_i(T_0, \tilde x, \a^{*,-i}, \a_i) \ge J_i(T_0, \tilde x, \a^*),\q \mbox{for all $i$, all $\tilde x\in \dbS_{T_0}$,  and all $\a_i \in \cA_{state, i}$}.
\eea
 We emphasize that the $\hat \a_i$ constructed in \reff{a1} is not in $\cA_{state,i}$, even when the $\a^*$ and $\a_i$ there are state dependent, so a more careful argument is required. We shall prove \reff{state-induction} by backward induction on $T_0$.

First, if $T_0 = T-1$, then the counterpart of \reff{a1} becomes: for any fixed $\tilde x\in \dbS_{T_0}$,
\beaa
\hat \a_i (s, \hat x) :=  \a_i  (s, \tilde x)\1_{\{s = T_0\}\cap \{\hat x = \tilde x\}} + \a^*_i(s, \hat x)\1_{\{s< T_0\}\cup\{\hat x \neq \tilde x\}},
\eeaa
which is in $\cA_{state, i}$.  Then \reff{state-induction} follows from the same arguments in Theorem \ref{thm-DPP-discrete} Step 1.

Assume \reff{state-induction}  holds true for $T_0+1$. Now for $T_0$, note that
\bea
\label{state-est1}
&&J_i(T_0, \tilde x, \a^{*,-i}, \a_i) \nonumber\\
&&= f_i(T_0, \tilde x, \a_i(T_0, \tilde x)) + \sum_{\hat x\in \dbS_{T_0+1}} q(T_0, \tilde x, (\a^{*,-i}, \a_i)(T_0, \tilde x), \hat x) J_i(T_0+1, \hat x, \a^{*,-i}, \a_i) \nonumber\\
&&\ge f_i(T_0, \tilde x, \a_i(T_0, \tilde x)) + \sum_{\hat x\in \dbS_{T_0+1}} q(T_0, \tilde x, (\a^{*,-i}, \a_i)(T_0, \tilde x), \hat x) J_i(T_0+1, \hat x, \a^*),
\eea
where the last inequality is due to the induction assumption. Fix $\tilde x\in \dbS_{T_0}$ and define
\beaa
\hat \a_i (s, \hat x) :=  \a_i  (s, \hat x)\1_{\{s = T_0\}\cap \{\hat x = \tilde x\}} + \a^*_i(s, \hat x)\1_{\{s\neq T_0\}\cup\{\hat x \neq \tilde x\}},
\eeaa
which is again state dependent. Then, denoting $\dbP^{t,x,\a}$ in the obvious way, 
\beaa
&&0 \le J_i(t, x, \a^{*,-i}, \hat \a_i) -  J_i(t,x, \a^*) = \dbP^{t,x,  \a^*}(X_{T_0} = \tilde x) \times\\
&&\qquad\q \Big[f_i(T_0, \tilde x, \a_i(T_0, \tilde x)) + \sum_{\hat x\in \dbS_{T_0+1}} q(T_0, \tilde x, (\a^{*,-i}, \a_i)(T_0, \tilde x), \hat x) J_i(T_0+1, \hat x, \a^*) \\
&&\qquad\qq- f_i(T_0, \tilde x, \a^*_i(T_0, \tilde x)) - \sum_{\hat x\in \dbS_{T_0+1}} q(T_0, \tilde x, \a^*(T_0, \tilde x), \hat x) J_i(T_0+1, \hat x, \a^*) \Big].
\eeaa
Note that $ \dbP^{t, x, \hat \a^*}(X_{T_0} = \tilde x) >0$. Then, together with \reff{state-est1}, the above implies 
\beaa
J_i(T_0, \tilde x, \a^{*,-i}, \a_i) &\ge& f_i(T_0, \tilde x, \a^*_i(T_0, \tilde x)) - \sum_{\hat x\in \dbS_{T_0+1}} q(T_0, \tilde x, \a^*(T_0, \tilde x), \hat x) J_i(T_0+1, \hat x, \a^*)\\
&=& J_i(T_0, \tilde x, \a^*).
\eeaa
This proves \reff{state-induction}, hence \reff{DPP-state3}.
 \qed

We now construct a counterexample such that the inclusion in \reff{DPP-state3} is strict. This is again due to the non-uniqueness of equilibriums.
\begin{eg}
\label{eg-state} 
 Let $T=4$, $N=2$, $A_1=A_2 = \{0,1\}$, and $\dbS^\dbT$ takes values as in Figure \ref{fig:state}.
 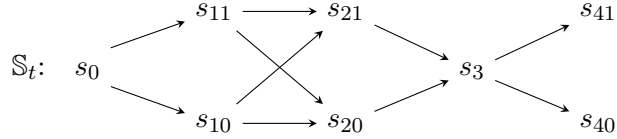
\begin{figure}[h]
\centering
\begin{tikzpicture}[>=stealth]

 \matrix (tree) [%
matrix of nodes,
minimum size=0.5cm,
column sep=1cm,
row sep=0.2cm,
]
{
                              &  $s_{11}$   & $s_{21}$ &            & $s_{41}$  \\
$\dbS_t$:\q $s_0$&                    &                & $s_3$ &                 \\
                             &  $s_{10}$    & $s_{20}$ &           & $s_{40}$ \\                  
};
\draw[->] (tree-2-1) -- (tree-1-2);
\draw[->] (tree-2-1) -- (tree-3-2);
\draw[->] (tree-1-2) -- (tree-1-3);
\draw[->] (tree-1-2) -- (tree-3-3);
\draw[->] (tree-3-2) -- (tree-1-3);
\draw[->] (tree-3-2) -- (tree-3-3);
\draw[->] (tree-1-3) -- (tree-2-4);
\draw[->] (tree-3-3) -- (tree-2-4);
\draw[->] (tree-2-4) -- (tree-1-5);
\draw[->] (tree-2-4) -- (tree-3-5);
\end{tikzpicture}
\caption{\label{fig:state} States for Example~\ref{eg-state}}
\end{figure}

\no We shall construct an equilibrium whose value is in $\tilde \dbV_{state}(0, s_0) \backslash \dbV_{state}(0, s_0)$. Set 
\beaa
T_0=1,\q  q(0,\cd) = {1\over 2},\q f(0,\cd) = 0.
\eeaa
Given a desired $\psi$, for any $\a\in \cA_{state}$, clearly  $J(1, \psi; 0, s_0, \a) = {1\over 2} [\psi(s_{10})+ \psi(s_{11})]$, and thus
\bea
\label{tildeVstate}
\tilde \dbV_{state}(0, s_0) =  \Big\{{1\over 2} [\psi(s_{10})+ \psi(s_{11})]: \mbox{for all $\psi$ s.t. $\psi(s_{1i}) \in \dbV_{state}(1, s_{1i})$, $i=0,1$}\Big\}.
\eea
Note that $\dbV_{state}(1, s_{10})$ and $\dbV_{state}(1, s_{11})$ are two different three-period games. Let the ($3$-period) subgames at branch $X_1=s_{10}$ and at branch $X_1=s_{11}$ exactly as in Example \ref{eg-path}. Since we consider only $\a\in \cA_{state}$,  by \reff{pathJ0} we have
\beaa
J(1, s_{1i}, \a) = \bar J(\a(3, s_3)), \q i=0,1.
\eeaa
Then, by Example \ref{eg-path}, 
\beaa
\dbV_{state}(1, s_{10}) = \dbV_{state}(1, s_{11}) = \big\{(0, {1\over 4}), ({1\over 4}, 0)\big\},
\eeaa
 with corresponding equilibriums $\a(3, s_3) = (0,0)$ and $\a(3,s_3) = (1,1)$ (the other values of $\a(t, \bx)$ is irrelevant, or say, can be arbitrary). Then, by \reff{tildeVstate},
\beaa
\tilde \dbV_{state}(0, s_0) = \big\{(0, {1\over 4}), ({1\over 4}, 0), ({1\over 8}, {1\over 8})\big\}.
\eeaa

On the other hand, since $q(0, \cd)={1\over 2}$ and $f(0,\cd) = 0$, for any $\a\in \cA_{state}$ we have
\beaa
J(0, s_0, \a) = {1\over 2}[J(1, s_{10}, \a) + J(1, s_{11}, \a)] = \bar J(\a(3, s_3)).
\eeaa
So $\dbV_{state}(0, s_0)= \{(0, {1\over 4}), ({1\over 4}, 0)\}$, therefore, $({1\over 8}, {1\over 8})\in \tilde \dbV_{state}(0, s_0) \backslash \dbV_{state}(0, s_0)$.
\end{eg}

\subsection{Pareto equilibriums}\label{sect-Pareto}

For $y, \tilde y \in \dbR^N$. We say $y \le \tilde y$ if $y_i \le \tilde y_i$ for $i=1,\cds, N$, and $y < \tilde y$ if we assume further that $y_i < \tilde y_i$ for some $i$.  As we saw in Remark \ref{rem-discrete-comparison} \ref{rem-discrete-comparison-2}, for a non-zero sum game typically the comparison principle fails in the sense: for equilibriums $\a^*,\tilde \a^*$ for games $J,\tilde J$, respectively,
\beaa
J(\a) \le \tilde J(\a) ~\mbox{for all $\a$, but}~ J(\a^*) > \tilde J(\tilde \a^*).
\eeaa


A consequence of the above property is that DPP would fail, in general, if one restricts to the so called Pareto equilibriums.
\begin{defn}
\label{defn-Pareto}
We say  $\a^*\in NE(t, \bx)$ is a Pareto equilibrium if there does not exist another equilibrium $\tilde \a\in NE(t,\bx)$ such that $J(t,\bx, \tilde \a) < J(t, \bx, \a^*)$.
\end{defn} 
Define
 \beaa
\dbV_{Pareto}(t,\bx) := \{J(t,\bx, \a^*): \mbox{for all Pareto equilibriums $\a^*\in NE(t,\bx)$}\}.
\eeaa
As the following example show, even the partial DPPs fail in general:
\bea
\label{DPP-Pareto}
\left.\ba{c}
\dbV_{Pareto}(t, \bx) \neq \Big\{J(T_0, \psi; t, \bx, \a^*): \mbox{for all $\psi$  and $\a^*$ such that } \\
\mbox{$\psi(\tilde\bx) \in \dbV_{Pareto}(T_0, \tilde \bx), \forall \tilde\bx\in \dbS^{\dbT}_{t, \bx}$, and $\a^*$ is a Pareto   equilibrium  at $(T_0, \psi; t, \bx)$}\Big\}.
\ea\right.
\eea

\begin{eg}
\label{eg-Pareto}  
As usual let $\tilde \dbV_{Pareto}(t, \bx)$ denote the right side of \reff{DPP-Pareto}. Let $T=2$, $N=2$, $A_1=A_2 = \{0,1\}$,  and $\dbS^\dbT$ takes values as in Figure \ref{fig:pareto}.
 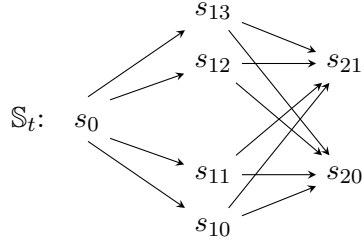
\begin{figure}[h]
\centering
\begin{tikzpicture}[>=stealth]

 \matrix (tree) [%
matrix of nodes,
minimum size=0.5cm,
column sep=1cm,
row sep=0.2cm,
]
{
                              &  $s_{13}$   &  		     \\
                              &  $s_{12}$   &  $s_{21}$   \\
$\dbS_t$:\q $s_0$&                    &                   \\
                             &  $s_{11}$    &  $s_{20}$   \\      
                             &  $s_{10}$    &                   \\              
};
\draw[->] (tree-3-1) -- (tree-1-2);
\draw[->] (tree-3-1) -- (tree-2-2);
\draw[->] (tree-3-1) -- (tree-4-2);
\draw[->] (tree-3-1) -- (tree-5-2);
\draw[->] (tree-1-2) -- (tree-2-3);
\draw[->] (tree-2-2) -- (tree-2-3);
\draw[->] (tree-4-2) -- (tree-2-3);
\draw[->] (tree-5-2) -- (tree-2-3);
\draw[->] (tree-1-2) -- (tree-4-3);
\draw[->] (tree-2-2) -- (tree-4-3);
\draw[->] (tree-4-2) -- (tree-4-3);
\draw[->] (tree-5-2) -- (tree-4-3);
\end{tikzpicture}
\caption{\label{fig:pareto} States for Example~\ref{eg-Pareto}}
\end{figure}

We first consider the subgame $\dbV(1, \bx)$. Set 
\beaa
 g(\bx)|_{\bx_2=s_{21}}=(0,0).
\eeaa
Let $f(1,x) := f(1,(s_0, x), a)$ (independent of $a$), $g(x) := g(s_0, x, s_{20})$, and $q(1,x, a):= q(1,\bx, a; s_{20})$ (independent of $\bx$) be as in Table \ref{tab:pareto1}. Then $J(1, x, a):= J(1, (s_0, x), a)$ are as in Table \ref{tab:pareto2}. This implies that
 \beaa
 \dbV(1, \bx) = \{\psi^*(\bx_1), \tilde \psi^*(\bx_1)\},\q  \dbV_{Pareto}(1, \bx) = \{\psi^*(\bx_1)\},
 \eeaa
 where $\psi^*$ and $\tilde \psi^*$ are given in Table \ref{tab:pareto3}.
 \begin{table}[h]
  \begin{center}
        \begin{tabular}{|c|c|c|c|c|} 
     \hline
      $x$ & $s_{10}$ & $s_{11}$ &  $s_{12}$ & $s_{13}$\\
      \hline
      $f(1,x)$ & $(1,1)$ & $(-4, 4)$ &   (4, -4)  & (1,1)\\
      \hline
      $g(x)$ & $(4,4)$ & $(20, 4)$ &    (4, 20)   & (12, 12)\\
      \hline
    \end{tabular}
     \qq
      \begin{tabular}{|c|c|c|} 
     \hline
      $q(1, a)$ & $a_2 = 0$ & $a_2 = 1$\\
      \hline
      $a_1=0$ & ${1\over 2}$ & ${3\over 4}$\\
      \hline
      $a_1=1$ & ${3\over 4}$ & ${1\over 4}$\\
      \hline
    \end{tabular}
    \end{center}
    \caption{\label{tab:pareto1} Cost and transition functions for Example~\ref{eg-Pareto}}
\end{table}

\begin{table}[h]
  \begin{center}
        \begin{tabular}{|c|c|c|} 
     \hline
      $J(1,s_{10}, a)$ & $a_2 = 0$ & $a_2 = 1$\\
      \hline
      $a_1=0$ & $(3,3)$ & $(4,4)$\\
      \hline
      $a_1=1$ & $(4,4)$ & $(2,2)$\\
      \hline
    \end{tabular}
    \qq
      \begin{tabular}{|c|c|c|} 
     \hline
      $J(1,s_{11}, a)$ & $a_2 = 0$ & $a_2 = 1$\\
      \hline
      $a_1=0$ & $(6,6)$ & $(11,7)$\\
      \hline
      $a_1=1$ & $(11,7)$ & $(1,5)$\\
      \hline
    \end{tabular}
    \bs\\
     \begin{tabular}{|c|c|c|} 
     \hline
      $J(1,s_{12}, a)$ & $a_2 = 0$ & $a_2 = 1$\\
      \hline
      $a_1=0$ & $(6,6)$ & $(7,11)$\\
      \hline
      $a_1=1$ & $(7,11)$ & $(5,1)$\\
      \hline
    \end{tabular}
    \qq
      \begin{tabular}{|c|c|c|} 
     \hline
      $J(1,s_{11}, a)$ & $a_2 = 0$ & $a_2 = 1$\\
      \hline
      $a_1=0$ & $(7,7)$ & $(10,10)$\\
      \hline
      $a_1=1$ & $(10,10)$ & $(4,4)$\\
      \hline
    \end{tabular}
  \end{center}
  \caption{\label{tab:pareto2} Cost matrices for Example~\ref{eg-Pareto}}
\end{table}
 
  \begin{table}[h]
  \begin{center}
        \begin{tabular}{|c|c|c|c|c|} 
     \hline
      $x$ & $s_{10}$ & $s_{11}$ &  $s_{12}$ & $s_{13}$\\
      \hline
      $\psi^*(x)$ & $(2,2)$ & $(1, 5)$ &   (5, 1)  & (4,4)\\
      \hline
      $\tilde \psi^*(x)$ & $(3,3)$ & $(6, 6)$ &    (6, 6)   & (7, 7)\\
      \hline
    \end{tabular}
    \end{center}
      \caption{\label{tab:pareto3} Values of the game in Example~\ref{eg-Pareto} at time $1$}
\end{table}

We now consider $J(1, \psi; 0, s_0, \a)$ for $\psi = \psi^*, \tilde \psi^*$.  Fix some $\e>0$ small enough. Set 
\beaa
f(0,\cd) = (0,0),\q  q(0, s_0, a; s_{1j}) = 1-3\e ~\mbox{if $j=I(a)$}\q \mbox{and}\q q(0, s_0, a; s_{1j}) = \e ~\mbox{if $j\neq I(a)$},
\eeaa
where
\beaa
I(0,0) =0,\q I(1,0) = 1,\q I(0,1) = 2, \q I(1,1) = 3.
\eeaa
and all other $q(0, s_0, a; x) = \e$. Then
\beaa
J(1,\psi; 0, s_0, a) = \sum_{j=0}^3 q(0, s_0, a; s_{1j}) \psi(s_{1j}) = \psi(s_{1 I(a)}) + O(\e).
\eeaa
That is, $J(1, \psi; 0, s_0, a)$ is approximately equal to $\psi(s_{1 I(a)})$ and, when $\e$ is small enough, the two subgames have the same equilibrium.  In particular, recall the $J$ and $\tilde J$ in Example \ref{eg-Pareto}, we see that  
\beaa
J(1, \psi^*; 0, s_0, a) = J(a) + O(\e),\q J(1, \tilde \psi^*; 0, s_0, a) = \tilde J(a) + O(\e).
\eeaa
Then, by Theorem \ref{thm-DPP-discrete}, 
\beaa
\tilde \dbV_{Pareto}(0, s_0) = \{(4,4) + O(\e)\},\q \dbV(0, s_0) =\{(3,3)+O(\e), (4,4)+O(\e)\},
\eeaa 
and thus $\dbV_{Pareto}(0,s_0) =  \{(3,3)+O(\e)\}$. This implies that $\tilde \dbV_{Pareto}(0, s_0)$ and  $\dbV_{Pareto}(0,s_0)$ do not include each other, namely partial DPP fails in both directions.
\end{eg}

\begin{rem}
\label{rem-Pareto} 
{\rm We emphasize that in Definition \ref{defn-Pareto} a Pareto equilibrium $\a^*$ is only compared to other equilibriums. In general it is possible that there exists another control $\a\in \cA$ (not an equilibrium) such that $J(t,\bx, \a) < J(t, \bx, \a^*)$, see Remark \ref{rem-discrete-comparison} \ref{rem-discrete-comparison-1}. We may call an equilibrium $\a^* \in \cA$ a strong Pareto equilibrium if there is no such  control $\a\in \cA$. Denote
\beaa
\dbV^{strong}_{Pareto}(t,\bx) := \{J(t,\bx, \a^*): \mbox{for all strong Pareto equilibriums $\a^*$}\}.
\eeaa
In general DPP fails for $\dbV^{strong}_{Pareto}$ too. 
\qed}
\end{rem}

\subsection{Optimal equilibriums }
\label{sect-optimal}
We now fix $x_0 \in \dbS_0$ and consider $\dbV(0,x_0)$. In practice it is  important to determine which equilibrium to implement. For this purpose we introduce a central planner, and assume the central planner is interested in minimizing:
\bea
\label{V0-discrete}
V_0 := \inf_{y \in \dbV(0, x_0)} \sum_{i=1}^N \l_i y_i = \inf\big\{\sum_{i=1}^N \l_i J_i(0, x_0, \a^*):   \a^*\in NE(0, x_0)\big\}.
\eea
where $\l_i\ge 0$ with $\sum_{i=1}^N\l_i = 1$. Such problems are natural, say, for social welfares. By Proposition \ref{prop-compact}, the problem \reff{V0-discrete} has an optimizer $y^*\in \dbV(0, x_0)$, and correspondingly there exists $\a^*\in NE(0, x_0)$. Note that, when $\l_i>0$ for all $i$, such $\a^*$ is automatically a Pareto equilibrium. We remark that in general neither $y^*$ nor $\a^*$ is unique, however,  the central planer is indifferent to them and thus can pick an arbitrary one. More importantly, in practice it is quite easy to implement such an equilibrium, as we explain below.

\begin{rem}
\label{rem-implementation} 
{\rm 
\begin{enumerate} 
\item Assume the central planner picks an optimal equilibrium $\a^*$ and recommend it to the players. As long as each player believes the others would follow the recommended one, it is in his/her best interest to follow the same $\a^*$ since it is an equilibrium.  Moreover, since $\a^*$ is a  Pareto optimal one (assuming $\l_i>0$ for all $i$), the players  are unlikely to make a collective decision to choose a different equilibrium. 

\item The problem is quite different from a ``dictatorship'' scenario, where the dictator wants to minimize  
\beaa
\tilde V_0 := \inf_{\a\in \cA} \sum_{i=1}^N \l_i J_i(0,x_0, \a).
\eeaa
Assume the problem $\tilde V_0$ has an optimal argument $\tilde \a^*$ and the dictator forces the players to follow it. However, since $\tilde \a^*$ is (in general) not an equilibrium, the individual players have no incentive to follow it even if they believe the others would do so. Consequently, the dictator has to use regulation/penalty (or other means) to force them to implement this strategy, which adds to the social cost.
\qed
\end{enumerate}
}
\end{rem}

\begin{rem}
\label{rem-moving}
{\rm
Since DPP fails for the Pareto equilibriums as detailed in Section~\ref{sect-Pareto}, the dynamic version of \reff{V0-discrete} will generally be time inconsistent.  {\color{black}In particular, this implies that there need not, and typically will not, exist a moving scalarization  (a moving objective parameterized by an adapted process $\l$) as in Feinstein \& Rudloff \cite{FR} so that $\a^*$ is a consistent equilibrium for this problem.}  Time inconsistency, therefore, implies that though a central planner may dictate a socially beneficial equilibrium at time 0, at some time $t$ this may no longer be an optimal equilibrium for the subgame over $[t, T]$.
\qed}
\end{rem}

\section{The continuous time model}
\label{sec:cont}
\setcounter{equation}{0}
{\color{black}In this section we extend our results to a continuous time setting. We shall consider a diffusion model with drift controls only. In this case all the involved probability measures are equivalent. The case with volatility controls may require new insights, especially in light of Remark \ref{rem-nondegenerate}, and is left for future research.
 }
 
\subsection{The nonzero sum game}
Let $[0, T]$ be the time horizon, $(\O, \dbF=\{\cF_t\}_{0\le t\le T}, \dbP_0)$ a filtered probability space, $B$ a $d$-dimensional $\dbP_0$-Brownian motion. Consider a game with $N$ players. Let $A = A_1\times \cds\times A_N$ be a convex domain in a Euclidean space,  and $\cA = \cA_1\times \cds\times \cA_N$ the set of $\dbF$-progressively measurable $A$-valued processes. The data of the game satisfy the following basic properties, where the boundedness assumption is mainly for simplicity. 

\begin{assum}
\label{assum-bound}
$(b, f): [0, T]\times \O\times A\to \dbR^d\times \dbR^N$ are  $\dbF$-progressively measurable and bounded; and  $\xi: \O\to \dbR^N$ is $\cF_T$-measurable and bounded.
\end{assum}
As usual we omit the variable $\o$ in $b, f, \xi$.  For each    $\a\in \cA$, define
\beaa
{d\dbP^{\a}\over d\dbP_0} := M^{\a}_{T}  := \exp\Big(\int_0^T b(s, \a_s)\cd dB_s - {1\over 2} \int_0^T |b(s,\a_s)|^2ds\Big).
\eeaa
At time $t$, each player has the value defined through the conditional expectation:
\beaa
J_i(t,\a) := \dbE^{\dbP^{\a}}_t\Big[\xi_i + \int_t^T f_i(s, \a^i_s)ds\Big],\q i=1,\cds, N.
\eeaa
We remark that we may replace the above expectation with some nonlinear operator through BSDEs, see Remark \ref{rem-BSDE} (ii) below.
 We say $\a^*\in \cA$ is a Nash equilibrium at $t$ if
\beaa
J_i(t, \a^*) \le J_i( t, \a^{*, -i}, \a^i),\q\mbox{$\dbP_0$-a.s.\ for all $i$ and all $\a^i \in \cA_i$},
\eeaa
and we introduce the set value:
\beaa
\cV_t :=\big\{J(t,\a^*): \mbox{for all Nash equilibrium $\a^*$ at $t$}\big\}.
\eeaa
We remark that the elements of $\cV_t$ are $\cF_t$-measurable $\dbR^N$-valued random variables, and we shall consider the localization in $\dbR^N$ in the next subsection.  

Given $T_0$ and $\eta \in \dbL^\infty(\cF_{T_0}; \dbR^N)$, denote
\beaa
J_i(T_0, \eta; t,\a) := \dbE^{\dbP^{\a}}_t\Big[\eta_i + \int_t^{T_0} f_i(s, \a^i_s)ds\Big],\q i=1,\cds, N,
\eeaa
and we define Nash equilibrium at $(T_0,\eta; t)$ in the obvious way. As such we then have the following DPP. {\color{black} We remark that this result does not even require the right continuity of $\dbF$.}

\begin{thm}
\label{thm-DPP-cont1}
Under Assumption \ref{assum-bound}, for any $0\le t< T_0\le T$ it holds
\bea
\label{DPP-cont1}
\cV_t := \big\{ J(T_0, \eta; t, \a^*): \mbox{for all $\eta\in \cV_{T_0}$ and all Nash equilibrium $\a^*$ at $(T_0, \eta; t)$}\big\}.
\eea
\end{thm}
\proof Let $\tilde \cV_t$ denote the right side of \reff{DPP-cont1}. First, for $J(t,   \a^*)\in \cV_t$, denote $\eta:= J(T_0, \a^*)$. For any $i$ and $\a^i\in \cA_i$, denote $\hat \a^i:= \a^i \1_{[0, T_0]} + \a^*\1_{(T_0, T]}$. It is clear that
\beaa
J_i(T_0, \eta; t, \a^{*,-i}, \a^i) = J_i(t, \a^{*,-i}, \hat \a^i) \ge J_i(t, \a^*) =  J_i(T_0, \eta; t, \a^*).
\eeaa
That is, $\a^*$ is a Nash equilibrium at $(T_0, \eta; t)$. Moreover, assume by contradiction that $\eta \notin \cV_{T_0}$, then there exist $i$ and $\a^i \in \cA_i$ such that $\dbP_0(E_i) >0$, where $E_i := \{J_i(T_0, \a^{*,-i}, \a^i) < J_i(T_0, \a^*)\}$. Denote $\hat \a^i := \a^* \1_{[0, T_0]} + \1_{(T_0, T]}[ \a^i \1_{E_i} + \a^*\1_{E_i^c}]$.  Then
\beaa
J_i(t, \a^{*,-i}, \hat \a^i) &=& \dbE^{\dbP^{\a^*}}_t\Big[\int_t^{T_0} f_i(s, \a^{*, i}_s) ds +  J_i(T_0, \a^{*,-i}, \a^i) \1_{E_i} + J_i(T_0, \a^*) \1_{E_i^c}\Big]\\
&<&  \dbE^{\dbP^{\a^*}}_t\Big[\int_t^{T_0} f_i(s, \a^{*, i}_s) ds +   J_i(T_0, \a^*) \Big] = J_i(t, \a^*).
\eeaa
This contradicts with the assumption that $\a^*$ is an equilibrium at $t$. Thus $\eta \in \cV_{T_0}$, and therefore $J(t, \a^*)\in \tilde \cV_t$.

Next, let $ J(T_0, \eta; t, \a^*) \in \tilde \cV_t$ with desired $(\eta, \a^*)$. Since $\eta \in \cV_{T_0}$,  $\eta = J(T_0, \tilde \a^*)$ for some equilibrium $\tilde \a^*$ at $T_0$. Denote $\hat \a^* := \a^* \1_{[0, T_0]} + \tilde \a^*\1_{(T_0, T]}$, and for any $i$ and $\a^i\in \cA_i$, denote $\hat \a^i :=   \a^i \1_{[0, T_0]} +\tilde \a^*\1_{(T_0, T]}$. Then
\beaa
&&J_i(t, \hat \a^{*, -i}, \a^i) - J_i(t, \hat \a^*) \\
&&= \big[J_i(t, \hat \a^{*, -i}, \a^i) -  J_i(t, \hat \a^{*,-i}, \hat \a^i) \big] +  \big[J_i(t, \hat \a^{*,-i}, \hat \a^i) - J_i(t, \hat \a^*)\big]\\
&& =\dbE^{\dbP^{\hat \a^{*, -i}, \a^i}}_t\big[ J_i(T_0, \tilde \a^{*, -i}, \a^i) -  J_i(T_0, \tilde \a^*) \big] + \big[J_i(T_0, \eta; t, \a^{*,-i}, \a^i) - J_i(T_0, \eta; t, \a^*)\big].
\eeaa
The second term above is nonnegative by the requirement of $\a^*$. Moreover,  note that $J_i(T_0, \tilde \a^{*, -i}, \a^i) \ge  J_i(T_0, \tilde \a^*)$, $\dbP_0$-a.s., and $\dbP^{\hat \a^{*, -i}, \a^i}$ is equivalent to $\dbP_0$, then $J_i(T_0, \tilde \a^{*, -i}, \a^i) \ge  J_i(T_0, \tilde \a^*)$, $\dbP^{\hat \a^{*, -i}, \a^i}$-a.s. This implies $J_i(t, \hat \a^{*, -i}, \a^i) \ge J_i(t, \hat \a^*)$. So $\a^*$ is an equilibrium at $t$, and thus $J(T_0, \eta; t, \a^*) = J(t, \hat \a^*) \in \cV_t$.
\qed

\subsection{The localization}
While Theorem \ref{thm-DPP-cont1} is quite simple, as mentioned $\cV_t$ is a set of random variables, rather than value sets in $\dbR^N$ as in Section \ref{sec:discrete}, which is not desirable in applications. In this subsection we localize the random variables in a pointwise sense. For this purpose it is more convenient to use the canonical space.

For the rest of this section, let  $\O := \{\o \in C([0, T]; \dbR^d): \o_0=0\}$ be the canonical space, $B$ the canonical process: $B_t(\o)=\o_t$,  $\dbP_0$ the Wiener measure, and $\dbF = \{\cF_t\}_{0\le t\le T}: = \dbF^B$ the $\dbP_0$-augmented filtration generated by $B$. Denote
\beaa
\|\o\| := \sup_{0\le t\le T} |\o_t|,\q {\bf d}((t, \o), (\tilde t, \tilde \o)) := \sqrt{|t-\tilde t|} + \|\o_{t\wedge \cd} - \tilde \o_{\tilde t\wedge \cd}\|.
\eeaa
Then $(\O, \|\cd\|)$ is a Polish space. For $t\in [0, T]$, $\o, \tilde \o\in \O$, and  $\xi\in \dbL^0(\cF_T)$, $\zeta\in \dbL^0(\dbF)$, denote 
\beaa
&(\o\otimes_t \tilde \o)_s := \o_s \1_{[0, t]}(s) + [\o_t + \tilde\o_{s-t}] \1_{[t, T]}(s),&\\
& \xi^{t,\o}(\tilde \o) := \xi(\o\otimes_t \tilde \o),\q \zeta^{t,\o}_s(\tilde \o) := \zeta_{t+s}(\o\otimes_t \tilde\o).&
\eeaa
Let $A, \cA, b,  f, \xi$ be as in the previous subsection. For $(t,\o) \in [0, T]\times \O$ and   $\a\in \cA$, define
\bea
\label{Ptoa}
\left.\ba{c}
\dis {d\dbP^{t,\o,\a}\over d\dbP_0} := M^{t,\o,\a}_{T-t}  := \exp\Big(\int_0^{T-t} b^{t,\o}(s, B_\cd, \a_s)\cd dB_s - {1\over 2} \int_0^{T-t} |b^{t,\o}(s, B_\cd, \a_s)|^2ds\Big);\\
\dis J_i(t,\o, \a) := \dbE^{\dbP^{t,\o,\a}}\Big[\xi^{t,\o}_i(B_\cd) + \int_0^{T-t} f^{t,\o}_i(s, B_\cd, \a^i_s)ds\Big],\q i=1,\cds, N.
\ea\right.
\eea
 We say $\a^*\in \cA$ is a Nash equilibrium at $(t,\o)$, denoted as $\a^*\in NE(t,\o)$, if
\beaa
J_i(t,\o, \a^*) \le J_i(t,\o, \a^{*, -i}, \a^i),\q\mbox{for all $i$ and all $\a^i \in \cA_i$},
\eeaa
and we introduce the set value:
\beaa
\dbV_0(t, \o) :=\big\{J(t,\o, \a^*): \a^*\in NE(t,\o)\big\} \subset \dbR^N.
\eeaa

Intuitively, $\eta\in \cV_t$ means $\eta(\o) \in \dbV_0(t,\o)$ for $\dbP_0$-a.e. $\o$. {\color{black}This is indeed true in the setting of Section \ref{sec:discrete} if we introduce the corresponding $\cV_t$. However, in the continuous time model we encounter some serious measurability issue.} Since the state space $\O$ is uncountable, the measurability or even certain regularity of the set value will be required.  Note that  $\cA$ is typically not compact, so the arguments in  Proposition \ref{prop-compact} do not work here. In fact, in this case neither the (Borel or analytic) measurability of the set $\dbV_0(t,\o)\subset \dbR^N$ for fixed $(t,\o)$ nor the $\dbF$-progressive measurability of the mapping $(t,\o)\to \dbV_0(t,\o)$ is clear to us. 
 To get around of this difficulty we  relax the equilibriums to approximating ones, which are usually sufficient in practice.

\begin{defn}
\label{defn-NEe}
We say $\a^\e\in \cA$ is an $\e$-equilibrium at $(t,\o)$, denoted as $\a^\e\in NE_\e(t,\o)$, if 
\beaa
J_i(t,\o, \a^\e) \le J_i(t,\o, \a^{\e, -i}, \a^i) + \e,\q\mbox{for all $i$ and all $\a^i \in \cA_i$}.
\eeaa
\end{defn}
Denote $O_\e(y):= \{\tilde y\in \dbR^N: |\tilde y- y|< \e\} \subset \dbR^N$,and define
\beaa
 \dbV(t,\o) := \bigcap_{\e>0} \dbV_\e(t,\o)~\mbox{where}~ \dbV_\e(t,\o) := \Big\{ y\in O_\e(J(t,\o; \a^\e)):  \a^\e\in NE_\e(t,\o)\Big\}.
\eeaa
Clearly $\dbV_0(t,\o) \subset \dbV(t,\o)$. Moreover, we have the following simple but important properties.
{\color{black}
\begin{prop}
\label{prop-compact2}
Let Assumption \ref{assum-bound} hold.
\begin{enumerate}
\item $\dbV_\e(t,\o)$ is bounded and open;
\item\label{prop-compact2-2} For any $\e' < \e$, the closure {\rm cl}$(\dbV_{\e'}(t,\o)) \subset \dbV_{\e}(t,\o)$;
\item $\dbV(t,\o)$ is compact. Moreover, $\dbV(t,\o) \neq \emptyset$ whenever $NE_\e(t,\o) \neq \emptyset$ for all $\e>0$.
\end{enumerate}
\end{prop}
\proof 
(i) This result is obvious. 

(ii) One can easily see that cl$(\dbV_{\e'}(t,\o)) \subset \{y\in O_{\e-\e'}(\tilde y):\tilde y \in \dbV_{\e'}(t,\o)\} \subset \dbV_{\e}(t,\o)$. 


(iii) Since $\dbV_\e(t,\o)$ is bounded, the ${\rm cl}(\dbV_\e(t,\o))$ is compact. By \ref{prop-compact2-2} we see that $\dbV(t,\o) = \cap_{\e>0} {\rm cl}(\dbV_\e(t,\o))$ is also compact. Moreover, again since each ${\rm cl}(\dbV_\e(t,\o))$ is compact, we see that $\dbV(t,\o) \neq \emptyset$ whenever ${\rm cl}(\dbV_\e(t,\o))\neq \emptyset$ for all $\e>0$.
\qed

\begin{rem}
\label{rem-existence}
{\rm (i) It is obvious that cl$(\dbV_0(t, \o)) \subset \dbV(t,\o)$, however, the inclusion could be strict. Note that $\dbV_0(t,\o) \neq \emptyset$ if and only if the game has a true equilibrium, while $\dbV(t,\o) \neq \emptyset$  can occur even if no equilibrium exists. Such a relaxation could be useful for more general games where a true equilibrium may not exist, see e.g., Frei \& dos Reis \cite{DF},   Buckdahn, Cardaliaguet, \& Rainer \cite{BCR}, and Lin \cite{Lin} for some results in this direction (the latter two use strategies instead of closed-loop controls though). 

{\color{black}(ii) When we view a stochastic control problem as a game with one player and denote its (standard) value function as $v(t,\o)$, then we always have $\dbV(t,\o) = \{v(t,\o)\}$, but $\dbV_0(t,\o)$ could be empty. Similarly for a two person zero sum game, the standard value function corresponds to $\dbV$, not $\dbV_0$.
}
}
\end{rem}

For the rest of the properties, we impose the following regularities.

\begin{assum}
\label{assum-reg}
\begin{enumerate}
\item\label{assum-reg-1} 
$b, f$ are uniformly continuous in {\color{black}$(t,\o)$ under ${\bf d}$} and $\xi$ is uniformly continuous in $\o$ under $\|\cd\|$, with a common modulus of continuity function $\rho_0$.

\item\label{assum-reg-2} 
$b, f$ are uniformly continuous in $a$.
\end{enumerate}
\end{assum}

We then have the regularity and stability of $\dbV$ in the  spirit of Feinstein \cite{Feinstein}. However, we note that \cite{Feinstein} considers the set of equilibriums, while we consider the set of values. Given $D_n\subset \dbR^N$, we define the set valued limits as in Aubin \& Frankowska \cite{AF}: 
\beaa
\left.\ba{c}
\dis\liminf_{n\to\infty} D_n  =  \Big\{y \in \dbR^N: \lim_{n \to \infty} \inf_{y_n \in D_n} |y-y_n| = 0 \Big\} \\
\dis \limsup_{n\to\infty} D_n  =  \Big\{y \in \dbR^N: \liminf_{n \to \infty} \inf_{y_n \in D_n} |y-y_n| = 0 \Big\}.
\ea\right.
\eeaa
That is, the limit inferior (superior) denotes the set of $y\in \dbR^N$ such that there exist $y_n\in D_n$ (resp.\ subsequence) satisfying $\lim_{n\to \infty} y_n = y$. 

\begin{thm}
\label{thm-reg}
Let Assumptions \ref{assum-bound} and  \ref{assum-reg} \ref{assum-reg-1}
 hold.
\begin{enumerate}

\item\label{thm-reg-0} {\color{black} For any $\e_1 < \e_2$, there exists $\d>0$ such that }
\bea
\label{Vereg}
{\color{black}\dbV_{\e_1} (\tilde t, \tilde \o) \subset \dbV_{\e_2}(t,\o)\q\mbox{for all $(t,\o), (\tilde t, \tilde \o)$ satisfying}~ {\bf d}((t,\o),(\tilde t, \tilde \o)) \le \d.}
\eea 

\item\label{thm-reg-1}
If ${\bf d}((t_n, \o^n), (t, \o))\to 0$, then $\dis \dbV(t,\o) = \bigcap_{\e>0} \big[\liminf_{n\to\infty} \dbV_\e(t_n, \o^n)\big] = \bigcap_{\e>0} \big[\limsup_{n\to\infty} \dbV_\e(t_n, \o^n)\big]$.

\item\label{thm-reg-2}
Assume $(b^n, f^n, \xi^n)$ satisfy Assumption \ref{assum-reg}  uniformly and define  $\dbV^n_\e(t,\o)$  in the obvious way. If $(b^n, f^n, \xi^n) \to (b, f, \xi)$ uniformly, then 
\beaa
\dbV(t,\o) = \bigcap_{\e>0} \big[\liminf_{n\to\infty} \dbV^n_\e(t, \o)\big] = \bigcap_{\e>0} \big[\limsup_{n\to\infty} \dbV^n_\e(t, \o)\big].
\eeaa
\end{enumerate}
\end{thm}
\proof 
(i) {\color{black}We first claim that there exists a modulus of continuity function $\rho$ such that
\bea
\label{Jreg}
  |J(t,\o, \a) - J(\tilde t, \tilde \o, \a)| \le \rho\big({\bf d}((t,\o), (\tilde t, \tilde \o))\big),\q\forall (t, \o), (\tilde t, \tilde \o), \forall \a.
\eea
Then, let ${\bf d}((t,\o),(\tilde t, \tilde \o)) \le \d$ and $y\in O_{\e_1}(J(\tilde t, \tilde \o, \a^{\e_1}))\subset \dbV_{\e_1} (\tilde t, \tilde \o)$ where $\a^{\e_1} \in NE_{\e_1}(\tilde t, \tilde \o)$. For any $i$ and $\a^i$, by \reff{Jreg} we have
\beaa
J_i(t,\o, \a^{\e_1}) &\le& J_i(\tilde t, \tilde \o, \a^{\e_1}) + \rho(\d) \le J_i(\tilde t, \tilde \o, \a^{\e_1,-i}, \a^i)+\e_1 + \rho(\d)\\
&\le& J_i(t,\o, \a^{\e_1,-i}, \a^i)+\e_1 +2 \rho(\d).
\eeaa
Choose $\d>0$ small enough such that $2 \rho(\d) \le \e_2-\e_1$, we see that $\a^{\e_1} \in NE_{\e_2}(t,\o)$. Moreover, by \reff{Jreg} again we have
\beaa
\Big|y- J_i(t,\o, \a^{\e_1})\Big|\le \Big|y- J_i(\tilde t, \tilde \o, \a^{\e_1})\Big| + \rho(\d) < \e_1 + \rho(\d) \le \e_2.
\eeaa
So $y\in \dbV_{\e_2}(t,\o)$, and hence \reff{Vereg} holds.

We next prove \reff{Jreg}.  By \reff{Ptoa} we have
\bea
\label{Jtoa2}
 J_i(t,\o, \a) = \dbE^{\dbP_0}\Big[M^{t,\o,\a}_{T-t}  \big[\xi^{t,\o}_i(B_\cd) + \int_0^{T-t} f^{t,\o}_i(s, B_\cd, \a^i_s)ds\big]\Big].
 \eea
 Similarly we have the representation for $J_i(\tilde t,\tilde \o, \a)$. Denote
 \beaa
OSC_\d(B):= \sup_{|s-\tilde s|\le \d} |B_s-B_{\tilde s}|,\q \rho'(\d):=  \dbE\Big[\rho_0^2\big( \d + OSC_\d(B)\big)\Big].
 \eeaa
 Assume without loss of generality that $t\le \tilde t$. Then, 
\beaa
&&\dis  \dbE\Big[\Big|\int_0^{T-t} f^{t,\o}_i(s, B_\cd, \a^i_s)ds - \int_0^{T-\tilde t} f^{\tilde t,\tilde \o}_i(s, B_\cd, \a^i_s)ds\Big|^2\Big]\\
&&\dis  \le C\dbE\Big[\int_0^{T-\tilde t} \big|f^{t,\o}_i(s, B_\cd, \a^i_s)- f^{\tilde t,\tilde \o}_i(s, B_\cd, \a^i_s)\big|^2ds +\big|\int_{T-\tilde t}^{T-t} f^{t,\o}_i(s, B_\cd, \a^i_s)ds\big|^2 \Big]\\
&&\dis\le C\dbE\Big[\int_0^{T-\tilde t} \rho_0^2\big({\bf d}\big((t+s, \o\otimes_t B_\cd), (\tilde t+s, \tilde \o\otimes_{\tilde t} B_\cd)\big)\big)ds \Big]+C\d^2\\
&&\dis\le C\dbE\Big[\int_0^{T-\tilde t} \rho_0^2\big({\bf d}\big((t, \o), (\tilde t, \tilde \o)\big)+ OSC_{\tilde t-t}(B)\big)ds \Big]+C\d^2\\
&&\dis\le C\dbE\Big[\int_0^{T-\tilde t} \rho_0^2\big(\d+ OSC_\d(B)\big)ds \Big]+C\d^2\le C\rho'(\d) + C\d^2.
\eeaa
Similarly, 
\beaa
 &&\dbE^{\dbP_0}\Big[|\xi^{t,\o}_i(B_\cd) - \xi^{\tilde t,\tilde \o}_i(B_\cd)|^2\Big]\le  \rho'(\d);\\
&&\dis  \dbE^{\dbP_0}\Big[\Big|\int_0^{T-t} b^{t,\o}(s, B_\cd, \a_s)dB_s - \int_0^{T-\tilde t} b^{\tilde t,\tilde \o}(s, B_\cd, \a_s)dB_s\Big|^2\Big]\\
&&\dis \q \le C\dbE^{\dbP_0}\Big[\int_0^{T-\tilde t} \big|b^{t,\o}- b^{\tilde t,\tilde \o}\big|^2(s, B_\cd, \a_s)ds +\int_{T-\tilde t}^{T-t} |b^{t,\o}(s, B_\cd, \a_s)|^2ds \Big]\\
&&\dis\q \le C\rho'(\d) + C\d;\\
&&\dis  \dbE^{\dbP_0}\Big[\Big|\int_0^{T-t} |b^{t,\o}(s, B_\cd, \a_s)|^2 ds - \int_0^{T-\tilde t} |b^{\tilde t,\tilde \o}(s, B_\cd, \a_s)|^2ds\Big|^2\Big]\\
&&\dis \q \le C\dbE^{\dbP_0}\Big[\int_0^{T-\tilde t} \big||b^{t,\o}|^2 - |b^{\tilde t,\tilde \o}|^2\big|^2 (s, B_\cd, \a_s)ds +\Big(\int_{T-\tilde t}^{T-t} |b^{t,\o}(s, B_\cd, \a_s)|^2ds\Big)^2 \Big]\\
&&\dis\q \le C\rho'(\d) + C\d^2.
\eeaa
We note that, since $b$ is bounded, for any $p\ge 1$, 
\bea
\label{Mest}
\sup_{\a\in\cA} \dbE^{\dbP_0}\Big[(M^\a_T)^p + (M^\a_T)^{-p}\Big] \le C_p< \infty.
\eea
Moreover, note that  $|e^x-e^{\tilde x}| \le [e^x + e^{\tilde x}] |x-\tilde x|$. Then
\beaa
&&\dis\dbE^{\dbP_0}\Big[ M^{t,\o,\a}_{T-t}\big[\big|\xi^{t,\o}_i(B_\cd) - \xi^{\tilde t,\tilde \o}_i(B_\cd)\big|\Big]\le  C\Big(\dbE^{\dbP_0}\Big[\big|\xi^{t,\o}_i(B_\cd) - \xi^{\tilde t,\tilde \o}_i(B_\cd)\big|^2 \Big]\Big)^{1\over 2}\le C\sqrt{\rho'(\d)};\\
&&\dis \dbE^{\dbP_0}\Big[ M^{t,\o,\a}_{T-t}\big|\int_0^{T-t} f^{t,\o}_i(s, B_\cd, \a^i_s)ds - \int_0^{T-\tilde t} f^{\tilde t,\tilde \o}_i(s, B_\cd, \a^i_s)ds\big|\big]\Big]\\
&&\dis\q \le C\Big(\dbE^{\dbP_0}\Big[\big|\int_0^{T-t} f^{t,\o}_i(s, B_\cd, \a^i_s)ds - \int_0^{T-\tilde t} f^{\tilde t,\tilde \o}_i(s, B_\cd, \a^i_s)ds\big|^2\Big]\Big)^{1\over 2}\\
&&\dis \q \le C\sqrt{\rho'(\d) + \d^2};\\
&&\dis \dbE^{\dbP_0}\Big[|M^{t,\o,\a}_{T-t}- M^{\tilde t,\tilde \o,\a}_{T-\tilde t}|\Big]\le \dbE^{\dbP_0}\Big[[M^{t,\o,\a}_{T-t}+ M^{\tilde t,\tilde \o,\a}_{T-\tilde t}]\times\\
&&\dis\q\big[\big|\int_0^{T-t} b^{t,\o}(s, B_\cd, \a_s)dB_s - \int_0^{T-\tilde t} b^{\tilde t,\tilde \o}(s, B_\cd, \a_s)dB_s\big| \\
&&\dis\q +  {1\over 2}\big|\int_0^{T-t} |b^{t,\o}(s, B_\cd, \a_s)|^2 ds-\int_0^{T-\tilde t} |b^{\tilde t,\tilde \o}(s, B_\cd, \a_s)|^2ds\big|\big] \Big]\\
&&\dis\q \le C\sqrt{\rho'(\d) + \d};\\
&&\dis |J(t,\o, \a) - J(\tilde t, \tilde \o, \a)| \le \dbE^{\dbP_0}\Big[M^{t,\o,\a}_{T-t}\big|\xi^{t,\o}_i(B_\cd) - \xi^{\tilde t,\tilde \o}_i(B_\cd)| \\
&&\dis\q +M^{t,\o,\a}_{T-t} \big|\int_0^{T-t} f^{t,\o}_i(s, B_\cd, \a^i_s)ds - \int_0^{T-\tilde t} f^{\tilde t,\tilde \o}_i(s, B_\cd, \a^i_s)ds\big|\\
&&\dis \q + |M^{t,\o,\a}_{T-t}- M^{\tilde t,\tilde \o,\a}_{T-\tilde t}| \big| \xi^{\tilde t,\tilde \o}_i(B_\cd) + \int_0^{T-\tilde t} f^{\tilde t,\tilde \o}_i(s, B_\cd, \a^i_s)ds\big|\Big]\\
&&\q\le C\sqrt{\rho'(\d) + \d} =: \rho(\d).
\eeaa
Clearly $\rho'$ and hence $\rho$ are modulus of continuity functions, we thus obtain \reff{Jreg}.
}

(ii) Denote $\d_n := {\bf d}((t_n, \o^n), (t, \o)) \to 0$. For any $\e_1 < \e_2$, by \reff{Vereg} and its proof we have 
\bea
\label{inclusion}
\dbV_{\e_1}(t, \o) \subset \dbV_{\e_2}(t_n, \o^n),~ \dbV_{\e_1}(t_n, \o^n) \subset \dbV_{\e_2}(t, \o),\q \mbox{whenever}~ 2\rho(\d_n) \le \e_2 - \e_1.
\eea
Now fix $\e_2$ and set $\rho(\d_n) \le {\e_2\over 4}$, we see that \reff{inclusion} holds for all $\e_1 \le {\e_2\over 2}$. This implies immediately that $\dbV(t,\o) \subset \dbV_{\e_2}(t_n, \o^n)$ and $ \bigcap_{\e_1>0}\big[\limsup_{n\to\infty} \dbV_{\e_1}(t, \o^n) \big] \subset \dbV_{\e_2}(t,\o)$. Now send $\e_2\to 0$ we have $ \bigcap_{\e_1>0}\big[\limsup_{n\to\infty} \dbV_{\e_1}(t, \o^n) \big] \subset   \dbV(t,\o) \subset \bigcap_{\e>0} \big[\liminf_{n\to\infty} \dbV_\e(t, \o^n)\big]$. Since the limit inferior is always contained in the limit superior, hence they are all equal.

(iii) Let $J^n$ be defined by \reff{Jtoa2}, but corresponding to $(b^n, f^n, \xi^n)$. It is clear that $\dis c_n:= \sup_{t,\o, \a} |[J^n-J](t,\o,\a)| \to 0$. Then the result follows similar arguments to \ref{thm-reg-1}.  
\qed

To study the measurability of the mapping $(t,\o)\mapsto \dbV(t,\o)$,  we introduce 
\beaa
\hat \dbV_\e(t,\o) := \bigcup_{\e'< \e} \dbV_{\e'} (t,\o).
\eeaa
It is clear that 
\beaa
\hat \dbV_\e(t,\o) \subset  \dbV_\e(t,\o) \subset \hat \dbV_{\tilde \e}(t,\o),\q\forall \e<\tilde \e,\q\mbox{hence}\q \dbV(t,\o) = \bigcap_{\e>0} \hat \dbV_\e(t,\o).
\eeaa
We then have the following result, which will be quite useful for the DPP below. 
{\color{black}
\begin{thm}
\label{thm-measurability}
Let Assumptions \ref{assum-bound} and  \ref{assum-reg} \ref{assum-reg-1} hold. 
 For any $\e>0$, any $\dbF$-stopping time $\t$, and any $\eta \in \dbL^0(\cF_\t)$,  the events $\{\o\in \O: \eta(\o) \in \hat \dbV_\e(\t(\o),\o)\}$ and $\{\o\in \O: \eta(\o) \in  \dbV(\t(\o),\o)\}$ are $\cF_\t$-measurable. 
\end{thm}
}
\proof First note that $\{\o: \eta(\o) \in \dbV(\t(\o),\o)\} = \bigcap_{n\ge 1} \{\o: \eta(\o) \in \hat\dbV_{1\over n}(\t(\o),\o)\}$, then the measurability for $\dbV$ clearly follows from the measurability for $\hat\dbV_{1\over n}$. We now prove the claimed measurability for $\hat \dbV_\e$  in three steps.

{\it Step 1.} We first show that, for any $t$ and any compact set $K\subset\subset \dbR^N$, the event $\{\o\in \O: K \subset \hat \dbV_\e(t,\o)\}$ is open (in terms of $\o$ under $\|\cd\|$), and thus is obviously $\cF_t$-measurable. Indeed, fix $\o$ such that $K \subset \hat \dbV_\e(t,\o) = \bigcup_{\e'<\e} \dbV_{\e'}(t,\o)$. Since $K\subset \dbR^N$ is compact and $\dbV_{\e'}(t,\o) \subset \dbR^N$ is open and increasing in $\e'$, there exists $\e_1<\e$ such that $K \subset \dbV_{\e_1}(t,\o)$. Now by \reff{Vereg} we see that there exists $\d>0$ such that  $K\subset  \dbV_{\e_1+\e\over 2}(t,\tilde\o)\subset \hat \dbV_\e(t, \tilde \o)$ whenever $\|\tilde\o_{t\wedge \cd}-\o_{t\wedge \cd}\|\le \d$.

{\it Step 2.} We next show the result when $\t\equiv t$ is a constant. Note that the set of closed balls in $\dbR^N$ with rational centers and rational radii is countable, numerated  as $\{K_i\}_{i\ge 1}$. Since $\hat \dbV_\e(t,\o)$ is open,  for $\eta\in  \dbL^0(\cF_t)$ one can easily verify that
 \beaa
 \{\o: \eta(\o) \in \hat\dbV_\e(t, \o)\}  = \bigcup_{i\ge 1} \Big(E_i \cap  \{\o: \eta(\o) \in K_i\}\Big) ,\mbox{where}~ E_i :=  \{\o: K_i\subset \hat \dbV_\e(t, \o)\}. 
 \eeaa
 Clearly $\{\eta\in K_i\}$ is $\cF_t$-measurable, and by Step 1 the events $E_i$ are also $\cF_t$-measurable, then so is the event $\{\o: \eta(\o) \in \hat\dbV_\e(t, \o)\}$.

 {\color{black}
 {\it Step 3.} We now consider stopping times $\t$. If $\t$ is discrete, namely taking only finitely many values: $t_1,\cds, t_n$, then
 \beaa
 \{\o: \eta(\o)\in \hat\dbV_\e(\t(\o),\o)\} = \bigcup_{i=1}^n\Big(  \{\o: \eta(\o)\in \hat\dbV_\e(t_i,\o)\} \cap \{\t=t_i\} \Big),
 \eeaa
 By Step 2,  $\{\o: \eta(\o)\in \hat\dbV_\e(t_i,\o)\}\in \cF_{t_i}$ for each $i$, then the above clearly implies  $\{\o: \eta(\o)\in \hat\dbV_\e(\t(\o),\o)\} \in \cF_\t$.

Now for general $\t$, there exist stopping times $\t_n\downarrow \t$ such that each $\t_n$ is discrete and $0\le \t_n - \t \le 2^{-n}T$. Choose an arbitrary sequence $\e_m \uparrow \e$. By \reff{Vereg} , for any $m$ we have
\beaa
&\dis\{\o: \eta(\o) \in \dbV_{\e_{m-1}}(\t(\o), \o)\} \subset \liminf_{n\to\infty} \{\o: \eta(\o) \in \hat\dbV_{\e_{m}}(\t_n(\o), \o)\}\\
&\dis \subset \limsup_{n\to\infty} \{\o: \eta(\o) \in \hat\dbV_{\e_{m}}(\t_n(\o), \o)\} \subset   \{\o: \eta(\o) \in \dbV_{\e_{m+1}}(\t(\o), \o)\}.
\eeaa
Send $m\to\infty$ and note that the first and the last terms above have the same limit, then the middle two terms have to converge to the same limit, namely
\bea
\label{hatVlimit}
\lim_{m, n\to\infty}  \{\o: \eta(\o) \in \hat\dbV_{\e_{m}}(\t_n(\o), \o)\}=\{\o: \eta(\o) \in \hat\dbV_\e(\t(\o), \o)\}. 
\eea
We already have $ \{\o: \eta(\o) \in \hat\dbV_{\e_{m}}(\t_n(\o), \o)\}\in \cF_{\t_n}$. Since $\dbF$ is right continuous and $\t_n \downarrow \t$, then $ \lim_{n\to\infty} \{\o: \eta(\o) \in \hat\dbV_{\e_{m}}(\t_n(\o), \o)\}\in \cF_\t$, and thus $\{\o: \eta(\o) \in \hat\dbV_\e(\t(\o), \o)\}\in \cF_\t$.
 \qed
 }

\subsection{Dynamic programming principle}
Given an $\dbF$-stopping time $\t$ and $\eta \in \dbL^\infty(\cF_\t; \dbR^N)$, one may consider the game on $[0, \t]$ with terminal condition $\eta$. In particular,
\bea
\label{JT0}
J_i(\t, \eta; t,\o,\a) := \dbE^{\dbP^{t,\o,\a}}\Big[\eta^{t,\o}_i + \int_0^{\t^{t,\o}-t} f^{t,\o}_i(s, B_\cd, \a^i_s)ds\Big],\q t\le \t(\o),~ i=1,\cds, N,
\eea
and we can define equilibrium and $\e$-equilibrium at $(\t, \eta; t,\o)$ in the obvious sense.  We now state our main result of this section, extending Theorem \ref{thm-DPP-discrete} to the continuous time model.
\begin{thm}
\label{thm-locDPP} Let Assumptions  \ref{assum-bound} and  \ref{assum-reg}  hold. For any $(t, \o)$ and {\color{black}any $\dbF$-stopping time $\t$ with $\t(\o)>t$}, we have
\bea
\label{locDPP}
\left.\ba{c}
\dis\dbV(t, \o) = \bigcap_{\e>0}\Big\{ y\in O_\e(J(\t, \eta; t,\o,  \a^\e)): ~\mbox{for all $\eta\in \dbL^\infty(\cF_\t;\dbR^N)$ and $\a^\e\in \cA$ } \\
\dis\mbox{such that $\a^\e\in NE_\e(\t, \eta; t,\o)$ and }~ \dbP_0(\eta^{t,\o} \notin \hat\dbV_\e(\t^{t,\o}, B^{t,\o}_\cd)) \le \e\Big\}.
\ea\right.
\eea
\end{thm}

To prove the theorem, we first need a lemma. 
{\color{black}
\begin{lem}
\label{lem-reg} 
\begin{enumerate}
\item\label{lem-reg-1} Let $\t$ be an $\dbF$-stopping time and $ \eta\in \dbL^\infty(\cF_\t;\dbR^N)$. For any $\d>0$, there exist a discrete $\dbF$-stopping time $\t_\d$ with $0\le \t_\d-\t\le \d$, and an $\eta_\d\in  \dbL^\infty(\cF_{\t_\d};\dbR^N)$ with the same bound as $\eta$,  such that
\bea
\label{etad}
\dbE^{\dbP_0}\big[ |\eta_\d - \eta|\big] \le \d,\q\mbox{and $\eta_\d$ is uniformly continuous in $\o$}.
\eea

\item\label{lem-reg-2} 
For any $\a\in \cA$ and $\d>0$, there exists discrete $\a^\d= \sum_{i=0}^{n-1} \a^\d_{t_i} \1_{[t_i, t_{i+1})}\in \cA$   such that
\beaa
\dbE^{\dbP_0}\Big[\int_0^T [|\a^\d_t - \a_t|\wedge 1] dt\Big] \le \d, \mbox{and each $\a^\d_{t_i}$ is uniformly continuous in $\o$}.
\eeaa
\end{enumerate}
\end{lem}
\proof  (i) The case $\t\equiv t$  follows the same approximations in Zhang \cite[Theorem 2.5.2]{Zhang} Steps 1-4, and in this case we actually have $\t_\d \equiv t$ as well. We now prove (i) for general stopping time $\t$. First, clearly there exists discrete $\t_\d$ such that $0\le \t_\d-\t\le \d$. Assume $\t_\d$ takes values $t_1,\cds, t_n$. Since the space $(\O, \|\cd\|)$ is Polish and thus $\dbP_0$ is tight, see e.g. Billingsley \cite{Billingsley},    then for each $i$ there exists a compact set{\footnote{More rigorously, we should first get an $\cF^B_{t_i}$-measurable set $E_i \subset \{\t_\d = t_i\}$ with $\dbP_0(\{\t_\d = t_i\} \backslash E_i) = 0$, and then apply \cite{Billingsley} to obtain the desired $K_i\subset\subset E_i$.}} $K_i\subset\subset \{\t_\d = t_i\}$ such that $K_i\in \cF_{t_i}$ and $\dbP_0(\{\t_\d = t_i\} \backslash K_i) < {\d\over 3C_0n}$, where $C_0$ is the bound of $\eta$. Then one may easily construct uniformly continuous functions $I_i\in \dbL^0(\cF_{t_i}; [0,1])$ such that $\dbE^{\dbP_0}[|I_i-\1_{K_i}|]\le  {\d\over 3C_0n}$. Next, note that $\eta\in \dbL^\infty(\cF_{\t_\d}; \dbR^N)$, then $\eta\1_{\{\t_\d=t_i\}}$ is $\cF_{t_i}$-measurable. Apply \reff{etad} for the deterministic time case, there exist $\eta_{i} \in  \dbL^\infty(\cF_{t_i};\dbR^N)$ with the same bound as $\eta$  such that
\beaa
\dbE^{\dbP_0}\big[ |\eta_{i} - \eta\1_{\{\t_\d=t_i\}}|\big] \le {\d\over 3n},\q\mbox{and $\eta_{i}$ is uniformly continuous in $\o$}.
\eeaa
Denote $\eta_\d:= \sum_{i=1}^n \eta_{i}I_i$. Then one can easily verify that  $\eta_\d$ is $\cF_{\t_\d}$ measurable, uniformly continuous, and
\beaa
&\dis\dbE^{\dbP_0}\big[ |\eta_\d - \eta|\big]\le \sum_{i=1}^n \dbE^{\dbP_0}\big[ |\eta_iI_i - \eta\1_{\{\t_\d=t_i\}}|\big]\\
&\dis\le \sum_{i=1}^n \dbE^{\dbP_0}\Big[ |\eta_i[I_i -\1_{K_i}]| + |\eta_i[\1_{K_i} -\1_{\{\t_\d=t_i\}}]| + |\eta_i - \eta\1_{\{\t_\d=t_i\}}|\1_{\{\t_\d=t_i\}}\Big] \\
&\dis\le \sum_{i=1}^n \Big[C_0{\d\over 3C_0 n} + C_0{\d\over 3C_0 n} + {\d\over 3n}\Big] = \d.
\eeaa
This proves \reff{etad} for the general stopping time $\t$.

(ii) First, denote $\a^R_t:= \a_t \1_{\{|\a_t|\le R\}}$. Then $\lim_{R\to\infty}  \dbE\Big[\int_0^T [|\a^R_t - \a_t| \wedge 1] dt\Big] =0.$
By otherwise choosing an $\a^R$, without loss of generality we assume $\a$ is bounded. Next, for each $n$, denote $t_i := {i\over n}T$, $i=0,\cds, n$.  Denote $\a^n_t := 0$, $t\in [t_0, t_1]$, and $\a^n_t := {n\over T} \int_{t_{i-1}}^{t_i} \a_s ds$, $t\in (t_i, t_{i+1}]$, $i=1,\cds, n-1$. Then $\dbE\Big[\int_0^T |\a^n_t - \a_t|  dt \Big] \le {\d\over 2}$ for $n$ large. Now fix such an $n$.  For each $\a^n_{t_i}$, by (i) we may construct uniformly continuous $\a^\d_{t_i}$ such that $\dbE[|\a^\d_{t_i} - \a^n_{t_i}|] \le {\d\over 2}$. Then clearly $\a^\d$ satisfies all the claimed properties.
\qed
}

\no{\bf Proof of Theorem \ref{thm-locDPP}.}  For notational simplicity, we assume $t=0$, then \reff{locDPP} becomes:
\bea
\label{locDPP0}
\left.\ba{c}
\dbV(0,0) =\tilde V(0,0):= \bigcap_{\e>0} \tilde V_\e(0,0)\q\mbox{where} \\
\begin{split}
\tilde V_\e(0,0) := \Big\{ y&\in  O_\e(J(\t, \eta; 0,0,  \a^\e)): ~\mbox{for all $\eta\in \dbL^\infty(\cF_\t;\dbR^N)$ and $\a^\e\in \cA$} \\
&\mbox{such that $\a^\e\in NE_\e(\t, \eta; 0,0)$ and }~ \dbP_0(\eta \notin \hat\dbV_\e(\t, B_\cd)) \le \e\Big\}.
\end{split}
\ea\right.
\eea

{\color{black}
{\it Step 1.} We first  prove the $\subset$ part.  Fix an arbitrary $y\in  \dbV(0,0)$. To show $y\in \tilde V(0,0)$, we fix an arbitrary $\e>0$. Let $\d>0$ be a small number which  will be specified later.

Since $y\in \dbV_\d(0,0)$,  there exists $\tilde\a^\d\in NE_\d(0,0)$ such that $|y -  J(0,0,\tilde\a^\d)|\le \d$.  For any $\d_1>0$, apply Lemma \ref{lem-reg} \ref{lem-reg-2}
  on $\tilde\a^\d$, there exists  $\a^\d = \sum_{i=0}^{n-1} \a^\d_{t_i}\1_{[t_i, t_{i+1})}\in \cA$ such that $\a^\d_{t_i}$ is uniformly continuous in $\o$ and $\dbE^{\dbP_0}\big[\int_0^T [|\tilde \a^\d_t - \a^\d_t|\wedge 1]dt\big]\le \d_1$. By Assumption \ref{assum-reg} \ref{assum-reg-2} 
  and \reff{Jtoa2}, for $\d_1$ small enough (depending on $\d$) we see that
  \bea
  \label{tildead}
  \a^\d\in NE_{2\d}(0,0)\q\mbox{and}\q  |y -  J(0,0,\a^\d)|\le 2\d.
  \eea
    Define 
  \bea
\label{eta}
  \eta^\d( \o) := J(\t(\o), \o, (\a^\d)^{\t(\o), \o}). 
    \eea
  By \reff{Jreg} and Assumption \ref{assum-reg} \ref{assum-reg-2} again, it is clear that $\eta^\d$ is $\cF_\t$-measurable.  
   Note that, for any $\a\in\cA$, $J(\t, \eta^\d;0,0,\a) = J(0,0, \tilde \a)$, where $\tilde \a:= \a \1_{[0, \t]} + \a^\d \1_{(\t, T]}$. Then \reff{tildead} implies $\a^\d \in NE_{2\d}(\t, \eta^\d, 0,0)$ and $|y - J(\t, \eta^\d; 0,0,\a^\d)|\le 2\d$. We shall always set $2\d\le \e$. Moreover,  set $\e_m\uparrow \e$ and $\t_n\downarrow \t$ be as in Theorem \ref{thm-measurability} Step 3. 
   We claim that, for any $m$,
   \bea
   \label{hatVlimit2}
   \limsup_{\d\to 0}\limsup_{n\to\infty}  \dbP_0\big(\{\o:\eta^\d(\o) \notin \hat \dbV_{\e_m}(\t_n(\o),\o)\}\big) =0.
   \eea
   Then, by \reff{hatVlimit} and noting that $\hat\dbV_\e$ is increasing in $\e$, we can easily see that 
   \beaa
  &\dis \limsup_{\d\to 0}\dbP_0\big(\{\o:\eta^\d(\o) \notin \hat \dbV_\e(\t(\o),\o)\}\big) =\limsup_{\d\to 0}\lim_{m,n\to\infty}\dbP_0\big(\{\o:\eta^\d(\o) \notin \hat \dbV_{\e_m}(\t_n(\o),\o)\}\big)\\
  &\dis \le \limsup_{\d\to 0}\limsup_{n\to\infty}\dbP_0\big(\{\o:\eta^\d(\o) \notin \hat \dbV_{\e_1}(\t_n(\o),\o)\}\big)=0. 
   \eeaa
This verifies all the requirements  in \reff{locDPP0} and thus $y\in\tilde \dbV_\e(0,0)$.

We now prove \reff{hatVlimit2} for $m=1$.   Since $\a^\d$ is uniformly continuous in $\o$, by  \reff{JT0} and Assumption \ref{assum-reg} \ref{assum-reg-2}, similar to \reff{Jreg} we have
 \bea
 \label{etanconv}
\lim_{n\to\infty} \dbE\big[|\eta^\d_n - \eta^\d|\big]=0,\q\mbox{where}\q  \eta^\d_n(\o) := J(\t_n(\o), \o, (\a^\d)^{\t_n(\o), \o}).
\eea
Note that
\beaa
&\dis \big\{ \o: \eta^\d(\o) \notin \hat \dbV_{\e_1}(\t_n(\o),\o)\big\} \subset \big\{ \o: \eta^\d(\o) \notin \dbV_{\e_2}(\t_n(\o),\o)\big\} \\
&\dis\subset  E^\d_n \cup  \big\{\o: |\eta^\d(\o) - \eta^\d_n(\o)|>\e_2\big\},
\eeaa
where, assuming $\t_n$ takes values $t_i$, $i=0,\cds, 2^n$, 
\bea
\label{Ei}
E^\d_n := \bigcup_{i=0}^{2^n} E_i,\q E_i :=\{\t_n = t_i\}\cap  \big\{ \o:  (\a^\d)^{t_i,\o} \notin NE_{\e_2}(t_i,\o)\big\}.
\eea
Then
\beaa
\dbP_0 \big(\big\{ \o: \eta^\d(\o) \notin \hat \dbV_{\e_1}(\t_n(\o),\o)\big\}\big) \le \dbP_0(E^\d_n) + {1\over\e_2}\dbE^{\dbP_0}\Big[|\eta^\d_n -\eta^\d(\o) |\Big].
\eeaa
By \reff{etanconv}, it suffices to show that
\bea
  \label{hatVlimit3}
   \limsup_{\d\to 0}\limsup_{n\to\infty}  \dbP_0(E^\d_n) =0.
   \eea

Now fix $\d, n$. Let $\d'>0$ be another small number which will be specified later. Note that $\O$ is separable, we may have a decomposition $E_i = \cup_{j\ge 1} E_{i,j}$ on $\cF_{t_i}$ such that, for some fixed $\o^{i,j}\in E_{i,j}$, $\sup_{0\le s\le t_i} |\o_s-\o^{i,j}_s|\le \d'$ for all $\o\in E_{i,j}$. Now for each $(i,j)$, since $(\a^\d)^{t_i,\o^{i,j}} \notin NE_{\e_2}(t_i,\o^{i,j})$, there exists $k=1,\cds, N$ and $\a^{i,j,k}\in \cA$ such that
\bea
\label{Jkij1}
(\eta^\d_n)_k(\o^{i,j}) = J_k(t_i,\o^{i,j}, (\a^\d)^{t_i,\o^{i,j}}) > J_k(t_i,\o^{i,j}, (\a^{\d,-k})^{t_i,\o^{i,j}}, \a^{i,j,k})  + \e_2.
\eea
Again by  \reff{JT0} and Assumption \ref{assum-reg} \ref{assum-reg-2}, and since $\a^\d$ is uniformly continuous in $\o$, then $\o\mapsto J_k(t_i,\o, (\a^\d)^{t_i,\o})$ and $\o\mapsto J_k(t_i,\o, (\a^{\d,-k})^{t_i,\o}, \a^{i,j,k})$ are uniformly continuous. Thus, for $\d'$ small enough,
\bea
\label{Jkij2}
J_k(t_i,\o, (\a^\d)^{t_i,\o}) > J_k(t_i,\o, (\a^{\d,-k})^{t_i,\o}, \a^{i,j,k})  + {\e_2\over 2},\q\forall \o\in E_{i,j}.
\eea
Denote $E^k_{i,j} := \{\o\in E_{i, j}: \mbox{\reff{Jkij1} holds}\}$. Then $E^\d_n = \bigcup_{k=1}^N E^k$, where $E^k:= \bigcup_{i=0}^{2^n} \bigcup_{j\ge 1} E^k_{ij}$. One can easily construct $\a^k\in\cA_k$ such that $(\a^k)_t^{t_i, \o} = \a^{i,j,k}_t$ for $(t,\o)\in [t_i, T]\times E^k_{i,j}$, and $\a^k_t= (\a^\d)^k_t$ for all other $(t, \o)$. Then by \reff{Jkij2} we have
\beaa
&\dis J_k(\t_n(\o),\o, (\a^\d)^{\t_n(\o),\o}) > J_k(\t_n(\o),\o, (\a^{\d,-k})^{\t_n(\o),\o}, (\a^k)^{\t_n(\o),\o})  + {\e_2\over 2},\q\forall \o\in E^k;\\
&\dis J_k(\t_n,\o, (\a^\d)^{\t_n(\o),\o}) = J_k(\t_n(\o),\o, (\a^{\d,-k})^{\t_n(\o),\o}, (\a^k)^{\t_n(\o),\o}),\q\forall \o\notin E^k.
\eeaa
Note that $\a^k_t = (\a^\d)^k_t$ for $t\le \t_n$. Then, since $\a^\d\in NE_{2\d}(0,0)$, 
\beaa
2\d &\ge& J_k(0,0, \a^\d) - J_k(0,0, \a^{\d, -k}, \a^k) \\
&=& \dbE^{\dbP^{0,0,\a^\d}}\Big[J_k(\t_n(\o), \o, (\a^\d)^{\t_n(\o),\o}) - J_k(\t_n(\o),\o, (\a^{\d,-k})^{\t_n(\o),\o}, (\a^k)^{\t_n(\o),\o})\Big]\\
&\ge& {\e_2\over 2} \dbP^{0,0,\a^\d}(E^k)={\e_2\over 2} \dbE^{\dbP_0} \Big[M^{0,0,\a^\d}_{\t_n}\1_{E^k}\Big].
\eeaa
Thus, by \reff{Mest},
\beaa
\dbP_0(E^k) = \dbE^{\dbP_0}\Big[ (M^{0,0,\a^\d}_{\t_n})^{-{1\over 2}} (M^{0,0,\a^\d}_{\t_n})^{{1\over 2}} \1_{E^k}\Big] \le C\Big(\dbE^{\dbP_0} \big[M^{0,0,\a^\d}_{\t_n}\1_{E^k}\big]\Big)^{1\over 2} \le C\sqrt{\d\over \e_2}.
\eeaa
Then $\dbP_0(E^\d_n) \le CN\sqrt{\d\over \e_2}$. This implies \reff{hatVlimit3} and hence \reff{hatVlimit2} immediately.

{\it Step 2.}  To see the opposite inclusion, we fix $y\in \tilde \dbV(0,0)$ and $\e>0$. Let $\d >0$ be a small number which will be specified later. Since $y\in \tilde V_\d(0,0)$, let $\eta, \a^\d$ be the corresponding terms in \reff{locDPP0} corresponding to $\d$. Moreover, set $\d_n \downarrow 0$ and let $(\t_n, \eta_n)$ be the approximations of $(\t, \eta)$ as in Lemma \ref{lem-reg} \ref{lem-reg-1} with error $\d_n$. Note that, for any $k=1,\cds, N$ and any $\a^k\in \cA_k$,
\bea
\label{Jktaun1}
&&\dis \Big|J_k(\t_n, \eta_n; 0,0, \a^{\d,-k}, \a^k) - J_k(\t, \eta; 0,0, \a^{\d,-k}, \a^k)\Big|\nonumber\\
&&\dis = \Big|\dbE^{\dbP_0}\Big[M^{ \a^{\d,-k}, \a^k}_{\t_n} \Big[[(\eta_n)_k - \eta_k]  + \int_\t^{\t_n}f_k(s, B, \a_s^{\d,-k}, \a^k_s) ds \Big]\Big| \\
&&\dis \le C\Big(\dbE^{\dbP_0}[|\eta_n - \eta|]\Big)^{1\over 2} + C2^{-n}\le \d,\nonumber
\eea
when $n$ is large enough. Thus
\bea
\label{Jktaun2}
\left.\ba{c}
\dis J_k(\t_n, \eta_n; 0,0, \a^{\d}) - J_k(\t_n, \eta_n; 0,0, \a^{\d,-k}, \a^k)\\
\dis \le J_k(\t, \eta; 0,0, \a^{\d})-J_k(\t, \eta; 0,0, \a^{\d,-k}, \a^k)+2\d \le 3\d,\q\forall \a^k\in \cA_k.
\ea\right.
\eea
That is, $\a^\d\in NE_{3\d}(\t_n, \eta_n; 0, 0)$ and $y\in O_{2\d}(J(\t_n, \eta_n; 0,0, \a^{\d}))$ for $n$ large enough.

Next, by \reff{hatVlimit} and noting that $\hat\dbV_\d$ is increasing in $\d$,  we have
\beaa
\limsup_{n\to\infty} \dbP_0(\eta \notin \hat\dbV_\d(\t_n, B_\cd)) \le \dbP_0(\eta \notin \hat\dbV_\d(\t, B_\cd))  \le \d.
\eeaa
Note that $\{\eta\in \hat \dbV_\d(\t_n, B)\}\cap \{|\eta_n - \eta|\le \d\} \subset \{\eta_n \in \hat \dbV_{2\d}(\t_n, B)\}$, then
\beaa
\limsup_{n\to\infty} \dbP_0(\eta_n \notin \hat\dbV_{2\d}(\t_n, B_\cd))\le \limsup_{n\to\infty} \Big[\dbP_0(\eta \notin \hat\dbV_{\d}(\t_n, B_\cd)) + {1\over \d} \dbE^{\dbP_0}[|\eta_n - \eta|]\Big]\le\d.
  \eeaa
 Thus, for $n$ large enough,
 \bea
\label{hatVlimit4}
 \dbP_0((E^\d_n)^c)\le 2\d,\q\mbox{where}\q E^\d_n := \{\eta_n \in \hat\dbV_{2\d}(\t_n, B_\cd)\}.
  \eea
  
Now fix $\d, n$. Denote  $E_i := E^\d_n\cap \{\t_n = t_i\}$. Let $\d'>0$ be another small number which will be specified later. Similar to Step 1 we have decomposition $E_i = \cup_{j\ge 1} E_{i,j}$ on $\cF_{t_i}$ such that,  for some fixed $\o^{i,j}\in E_{i,j}$, $\sup_{0\le s\le t_i} |\o_s-\o^{i,j}_s|\le \d'$ for all $\o\in E_{i,j}$. Now for each $(i, j)$, since $\eta_n(\o^{i,j})\in  \hat\dbV_{2\d}(t_i, \o^{i,j})\subset \dbV_{2\d}(t_i, \o^{i,j})$, there exists $\a^{i,j}\in \cA$ such that 
\beaa
\a^{i,j}\in NE_{2\d}(t_i, \o^{i,j})),\q |\eta_n(\o^{i,j})- J(t_i, \o^{i,j}, \a^{i,j})|\le 2\d.
\eeaa
 Since $\eta_n$ and $J(t_i, \o, \a^{i,j})$ are uniformly continuous in $\o$, for $\d'$ small enough we have 
\beaa
|\eta_n(\o)- \eta_n(\o^{i,j})|\le \d,\q \sup_\a |J(t_i, \o, \a)-J(t_i, \o^{i,j}, \a)|\le \d,\q\forall \o\in E_{i,j}.
\eeaa
Denote
\bea
\label{etantilde}
\tilde \eta_n(\o):= \sum_{i,j} \1_{E_{i,j}}(\o) J(t_i, \o, \a^{i,j})+  \1_{(E^\d_n)^c}(\o) J(\t_n(\o), \o, (\a^\d)^{\t_n(\o), \o}).
\eea
Then 
\beaa
&&\dis |\tilde \eta_n - \eta_n|\le  \sum_{i,j} \1_{E_{i,j}}(\o) |J(t_i, \o, \a^{i,j}) - \eta_n(\o)|+ C\1_{(E^\d_n)^c}\\
&&\dis \le  \sum_{i,j} \1_{E_{i,j}}(\o)\Big[|J(t_i, \o, \a^{i,j}) -J(t_i, \o^{i,j}, \a^{i,j})| +|J(t_i, \o^{i,j}, \a^{i,j})- \eta_n(\o^{i,j})|\\
&&\dis\q + |\eta_n(\o^{i,j})-\eta_n(\o)|\Big]+ C\1_{(E^\d_n)^c}\le 4\d + C\1_{(E^\d_n)^c}.
\eeaa
Similar to   \reff{Jktaun1} and  \reff{Jktaun2}, by \reff{hatVlimit4} one can easily show that 
\bea
\label{etantilde2}
\a^\d\in NE_{C\sqrt{\d}}(\t_n, \tilde \eta_n; 0,0)\q\mbox{and}\q y\in O_{C\sqrt{\d}}(J(\t_n, \tilde \eta_n; 0,0, \a^\d)).
\eea

We now define
\bea
\label{hatand}
\a^{\d,n}_t := \a^\d_t \1_{[0, \t_n]}(t) + \1_{(\t_n, T]}(t) \Big[\sum_{i,j}\1_{E_{i,j}} \a^{i,j}_{t-\t_n} + \a^\d_t\1_{(E^\d_n)^c}\Big]. 
\eea
Then $\tilde \eta_n(\o) = J(\t_n(\o), \o, (\a^{\d,n})^{\t_n(\o),\o})$ for all $\o\in\O$.  For $k=1,\cds, N$ and for $\a^k\in \cA_k$,  
\beaa
&&\dis J_k(0,0, \a^{\d,n}) =  J_k(\t_n, \tilde \eta_n; 0, 0, \a^\d) \le J_k(\t_n, \tilde \eta_n; 0, 0, \a^{\d,-k}, \a^k) + C\sqrt{\d}\\
&&\dis  =\dbE^{\dbP^{\a^{\d,-k}, \a^k}}\Big[\tilde \eta_n(\o) + \int_0^{\t_n} f_k(s, \o, \a^{\d,-k}_s,\a^k_s)ds\big]\Big]+ C\sqrt{\d}\\
&&\dis  \le\dbE^{\dbP^{\a^{\d,-k}, \a^k}}\Big[\sum_{i, j}\1_{E_{i,j}}J_k(t_i, \o, \a^{i,j}) + C\1_{(E^\d_n)^c}  + \int_0^{\t_n} f_k(s, \o, \a^{\d,-k}_s,\a^k_s)ds\big]\Big]+ C\sqrt{\d}\\
&&\dis  \le\dbE^{\dbP^{\a^{\d,-k}, \a^k}}\Big[\sum_{i, j}\1_{E_{i,j}}J_k(t_i, \o^{i,j}, \a^{i,j}) + \int_0^{\t_n} f_k(s, \o, \a^{\d,-k}_s,\a^k_s)ds\big]\Big]+ C\sqrt{\d}\\
&&\dis  \le\dbE^{\dbP^{\a^{\d,-k}, \a^k}}\Big[\sum_{i, j}\1_{E_{i,j}}J_k(t_i, \o^{i,j}, \a^{i,j,-k}, (\a^k)^{t_i, \o}) + \int_0^{\t_n} f_k(s, \o, \a^{\d,-k}_s,\a^k_s)ds\big]\Big]+ C\sqrt{\d}\\
&&\dis  \le\dbE^{\dbP^{\a^{\d,-k}, \a^k}}\Big[\sum_{i, j}\1_{E_{i,j}}J_k(t_i, \o, \a^{i,j,-k}, (\a^k)^{t_i, \o}) + \int_0^{\t_n} f_k(s, \o, \a^{\d,-k}_s,\a^k_s)ds\big]\Big]+ C\sqrt{\d}\\
&&\dis  \le\dbE^{\dbP^{\a^{\d,-k}, \a^k}}\Big[J_k(\t_n(\o), \o, (\a^{\d,n,-k},\a^k)^{\t_n(\o), \o}) + \int_0^{\t_n} f_k(s, \o, \a^{\d,-k}_s,\a^k_s)ds\big]\Big]+ C\sqrt{\d}\\
&&\dis = J_k(0, 0, \a^{\d,n,-k}, \a^k) +  C\sqrt{\d}.
\eeaa
That is, $\a^{\d, n}\in NE_{C\sqrt{\d}}(0,0)$, and $y\in O_{C\sqrt{\d}}(J(\t_n, \tilde \eta_n; 0,0, \a^\d)) =O_{C\sqrt{\d}}(J(0,0, \a^{\d,n}))$. Then $y\in \dbV_{C\sqrt{\d}}(0,0)$, and thus $y\in \dbV_\e(0,0)$ when $C\sqrt{\d}\le \e$.
\qed
}

\begin{rem}
\label{rem-cont-state}
{\rm In the state dependent setting, namely 
\bea
\label{cont-state}
b = b(t, \o_t, a),\q f=f(t, \o_t, a),\q \xi = g(\o_T),
\eea
as in Subsection \ref{sect-discrete-state} we can show that $\dbV(t,\o) = \dbV(t, \o_t)$ is also state dependent, but the DPP still involves path dependent $\eta$ and $\a^\e$.
\qed}
\end{rem}

\subsection{A duality result}
In this subsection we provide an alternative characterization for the  set value $\dbV(t,\o)$. The idea is similar to the level set or nodal set approach, see e.g. Barles, Soner, \& Souganidis \cite{BSS},  Ma \& Yong \cite{MY}, and Karnam, Ma, \& Zhang \cite{KMZ}. In particular, this method could be efficient for numerical purpose.

We first note that, for any $(t,\o)$ and $\a\in \cA$, $J(t,\o, \a) = Y^{t,\o,\a}_0$, where $(Y^{t,\o, \a}, Z^{t,\o, \a})$ is the solution to the following (linear) BSDE on $[0, T-t]$:
\bea
\label{BSDE}
\left.\ba{c}
\dis Y^{t,\o, \a,i}_s = \xi^{t,\o}_i(B) + \int_s^{T-t} f_i^{t,\o}(r, B, \a_r, Z^{ t,\o,\a,i}_r) dr - \int_s^{T-t} Z^{t,\o, \a, i}_r dB_r,\\
\dis \mbox{where}\q  f_i(t,\o, a, z_i) :=  f_i(t,\o, a_i) + b(t,\o, a) z_i.
\ea\right.
\eea 
For each $i$ and $a^{-i} = (a_1,\cds, a_{i-1}, a_{i+1}, \cds, a_N)$, denote
\beaa
\underline f_i(t,\o, a^{-i}, z_i) := \inf_{a_i\in A_i} f_i(t,\o, a^{-i}, a_i, z_i).
\eeaa
Since $b$ is bounded, $\underline f_i$ is uniformly Lipschitz continuous in $z_i$. {\color{black}Introduce the following multidimensional BSDE: $i=1,\cds, N$,
\bea
\label{BSDEunderline}
 \underline Y^{t,\o, \a,i}_s = \xi^{t,\o}_i(B) + \int_s^{T-t} \underline f_i^{t,\o}(r, B, \a^{-i}_r, \underline Z^{ t,\o,\a,i}_r) dr - \int_s^{T-t} \underline Z^{t,\o, \a, i}_r dB_r.
 \eea
 It is clear that, see e.g. El-Karoui \& Hamadene \cite{EH}, $\a^*\in NE(t,\o)$ if and only if 
 \bea
 \label{BSDENE}
  \underline f_i^{t,\o}(r, B, \a^{*,-i}_r, \underline Z^{ t,\o,\a^*,i}_r) =   f_i^{t,\o}(r, B,  \a^{*}_r, \underline Z^{ t,\o,\a^*,i}_r),~ a.s., 0\le r\le T-t, i=1,\cds, N.
  \eea

 Our main idea of the duality approach is to rewrite \reff{BSDEunderline} as a forward diffusion, by viewing the component $Z$ as a control. To be precise, }
fix $(t,\o, y)\in [0, T]\times \O\times \dbR^N$. For any   $\a\in \cA$ and $Z = (Z^1,\cds, Z^N)$, denote
\bea
\label{X}
Y^{t,\o, y, \a, Z,i}_s := y_i - \int_0^s \underline f_i^{t,\o}(r, B, \a^{-i}_r, Z^i_r) dr + \int_0^s Z^i_r dB_r.
\eea 
We then introduce an auxiliary control problem:
\bea
\label{W}
\left.\ba{c}
\dis W(t,\o, y) := \inf_{\a\in \cA, Z\in \dbL^2(\dbF, \dbP_0)}  \sum_{i=1}^N 
\dbE^{\dbP_0}\Big[ |\xi^{t,\o}_i(B) - Y^{t,\o, y,\a, Z, i}_{T-t}|^2 \\
\dis \qquad + \int_0^{T-t} [\D f^{t,\o}_i(s, B, \a_s, Z^i_s) ]^{3\over 2} ds\Big],\\
\dis \mbox{where}\q \D f_i(t, \o,  a, z_i):=  f_i(t, \o,  a, z_i)- \underline f_i(t,\o, a^{-i}, z_i).
\ea\right.
\eea
Here the power ${3\over 2}$ (between $1$ and $2$) for the $f$-term is  for some technical reasons on which we will elaborate later. {\color{black} By \reff{BSDEunderline}-\reff{BSDENE}, it is obvious that $W(t,\o, y) = 0$ for all $y\in \dbV_0(t,\o)$.}

Our main result of this subsection is that the set value agrees with the nodal set of $W$.
\begin{thm}
\label{thm-duality}
Let Assumptions \ref{assum-bound} and \ref{assum-reg}  hold. Then, for any $(t,\o)$, 
\beaa
\dbV(t,\o) = \dbN(t,\o) := \big\{y\in \dbR^N: W(t,\o, y) =0\big\}.
\eeaa
\end{thm}
\proof Without loss of generality, we assume $(t,\o) = (0,0)$, and for notational simplicity we may omit $(0,0)$ when there is no confusion, for example $J(\a) := J(0,0,\a)$. 

(i) We first show that $\dbN(0,0) \subset \dbV(0,0)$. Fix $y \in \dbN(0,0)$. For any $\e>0$, there exist $\a^\e$ and $Z^\e$ such that, denoting $Y^\e := Y^{y, \a^\e, Z^\e}$, 
\bea
\label{Ze}
\dbE^{\dbP_0}\Big[ |\xi_i - Y^{\e, i}_T|^2  + \int_0^T [ \D f_i(s, B, \a^{\e}_s, Z^{\e,i}_s) ]^{3\over 2} ds\Big] \le \e^2,~ i\ge 1.
\eea
 Let $(\tilde Y^\e, \tilde Z^\e)$ solve the following BSDE:
\beaa
\tilde Y^{\e,i}_t = \xi_i(B) + \int_t^T f_i(s, B, \a^{\e}_s, \tilde Z^{\e, i}_s) ds - \int_t^T \tilde Z^{\e, i}_s dB_s.
\eeaa
Note that
\bea
\label{BSDE2}
Y^{\e,i}_t = Y^{\e, i}_T + \int_t^T \underline f_i(s, B, \a^{\e, -i}_s, Z^{\e, i}_s) ds - \int_t^T Z^{\e, i}_s dB_s.
\eea
Then, denoting $\D Y^i:= \tilde Y^{\e,i} - Y^{\e,i}$ and $\D Z^i:= \tilde Z^{\e,i} - Z^{\e,i}$, we have
\beaa
\D Y^i_t = \xi_i(B) -  Y^{\e, i}_T + \int_t^T \D f_i(s, B, \a^{\e}_s, Z^{\e,i}_s)  ds+\int_t^T b(s, B,\a^\e_s)\D Z^i_sds - \int_t^T \D Z^{i}_s dB_r.
\eeaa
Thus, recalling \reff{Ptoa} for $M$, 
\beaa
\D Y^i_0 = \dbE^{\dbP_0}\Big[M^{\a^\e}_T\big[\xi_i(B) -  Y^{\e, i}_T + \int_0^T \D f_i(s, B, \a^{\e}_s, Z^{\e,i}_s) ds\big]\Big].
\eeaa
By \reff{Mest} and \reff{Ze} (in particular noting the power ${3\over 2}$ for the $f$-term is greater than $1$), it is clear that $|\tilde Y^{\e, i}_0 - Y^{\e, i}_0|\le C\e$. Moreover, let  $(\hat Y^\e, \hat Z^\e)$ solve the following BSDE:
\bea
\label{BSDE3}
\hat Y^{\e,i}_s = \xi_i(B)  + \int_s^T \underline f_i(s, B, \a^{\e, -i}_s, \hat Z^{\e, i}_s) dr - \int_s^T \hat Z^{\e, i}_r dB_r.
\eea
Compare \reff{BSDE2} and \reff{BSDE3}, it follows from \reff{Ze} again that $|\tilde Y^{\e, i}_0 -\hat Y^{\e, i}_0|\le C\e$, and thus $|\hat Y^{\e, i}_0 - Y^{\e, i}_0|\le C\e$.  

On the other hand, for any $\a^i$, applying the comparison principle on BSDEs \reff{BSDE} and \reff{BSDE3} we see that $J_i(\a^{\e,-i}, \a^i) \ge \hat Y^{\e,i}_0$. Then
\beaa
J_i(\a^\e) = Y^{\e,i}_0 \le \hat Y^{\e, i}_0 + C\e \le J_i(\a^{\e,-i}, \a^i) + C\e,
\eeaa
and thus $\a^\e\in NE_{C\e}(0, 0)$.  Recall $J(\a^\e) = Y^\e_0 = y$, then $y \in \dbV_{C\e}(0,0)$. Since $\e$ is arbitrary, we obtain $y\in \dbV(0,0)$. 

(ii) We next show that $\dbV(0,0) \subset \dbN(0,0)$.  Fix $y\in \dbV(0,0)$. For any $\e>0$, there exists $\a^\e\in NE_\e(0,0)$ such that $|y-J(\a^\e)| \le \e$. Recall that $J(\a^\e) = Y^{\a^\e}_0$, where $(Y^{\a^\e}, Z^{\a^\e})$ is defined by \reff{BSDE}. Let $(\hat Y^\e, \hat Z^\e)$ be defined by \reff{BSDE3}. For each $i$, there exists $\a^i$ such that 
\bea
\label{alphai}
f_i(r, B, \a^{\e,-i}_r, \a^i_r,\hat Z^{\e,i}_r ) \le \underline f_i(r, B, \a^{\e, -i},\hat Z^{\e, i}_r) + \e.
\eea
Let $(\check Y^{\e,i}, \check Z^{\e, i})$ solve the following BSDE:
 \bea
\label{BSDE4}
\check Y^{\e,i}_s = \xi_i(B) + \int_s^T f_i(r, B,  \a^{\e,-i}_r, \a^i_r, \check Z^{\e, i}_r) dr - \int_s^T \check Z^{\e, i}_r dB_r.
\eea
Compare BSDEs \reff{BSDE3} and \reff{BSDE4}, it follows from \reff{alphai} that $\check Y^{\e, i}_0 \le \hat Y^{\e,i}_0 + C\e$. Moreover, since $\a^\e\in NE_\e(0,0)$, then $Y^{\a^\e, i}_0 \le \check Y^{\e, i}_0 + \e\le \hat Y^{\e,i}_0 + C\e$. By the comparison principle of BSDEs we know that $Y^{\a^\e, i}_0 \ge \hat Y^{\e, i}_0$. Thus $|Y^{\a^\e, i}_0 - \hat Y^{\e,i}_0|\le C\e$. This, together with $|y- Y^{\a^\e}_0|\le \e$, implies that  $|y - \hat Y^{\e}_0|\le C\e$. 

Finally, note that 
\bea
\label{xiest}
 Y^{y, \a^{\e}, \hat Z^{\e}, i}_T - \xi_i(B)  = Y^{y, \a^{\e}, \hat Z^{\e}, i}_T - Y^{\hat Y^\e_0, \a^{\e}, \hat Z^{\e}, i}_T=  y_i-\hat Y^{\e,i}_0.
\eea
Moreover, note that $\underline f_i$ is uniformly Lipschitz in $z$. Then, denoting $\D Z^i := Z^{\a^\e, i} - \hat Z^{\e, i}$,
\beaa
&&\dis C\e \ge Y^{\a^\e, i}_0 -  \hat Y^{\e,i}_0 \\
&&\dis=  \int_0^T [f_i(s, B,  \a^\e_s, Z^{\a^\e, i}_s) - \underline f_i(s, B, \a^{\e, -i}_s, \hat Z^{\e,i}_s)] ds - \int_0^T \D Z^i_s dB_s\\
&&\dis =   \int_0^T \D f_i(s, B,  \a^\e_s, \hat Z^{\e, i}_s) ds  + \int_0^T \!\! b(s, B, \a^{\e}_s) \D  Z^{i}_s ds- \int_0^T\!\! \D Z^i_s dB_s.
\eeaa
This implies that
\bea
\label{ulfest}
\dbE^{\dbP_0}\Big[M^{\a^\e}_T \int_0^T \D f_i(s, B,  \a^\e_s, \hat Z^{\e, i}_s)ds\Big] \le C\e.
\eea
Since $\xi$ and $f$ are bounded, by standard BSDE estimates we have $\dbE^{\dbP_0}\Big[\int_0^T |\hat Z^{\e, i}_s|^2ds\Big] \le C$. Note further that
\beaa
0\le \D f_i(t,\o, a, z) \le C[1+|z|].
\eeaa 
One can easily derive from \reff{Mest} and \reff{ulfest} that (thanks to the fact ${3\over 2} < 2$)
\beaa
&&\dis\dbE^{\dbP_0}\Big[ \int_0^T [\D f_i(s, B,  \a^\e_s, \hat Z^{\e, i}_s)]^{3\over 2} ds\Big] \\
&&\dis\le C\dbE^{\dbP_0}\Big[ (M_T^{\a^\e})^{-{1\over 4}} (M_T^{\a^\e})^{1\over 4}\int_0^T [\D f_i(s, B,  \a^\e_s, \hat Z^{\e, i}_s)]^{1\over 4} ds  \int_0^T [1+|\hat Z^{\e, i}_s|^{5\over 4}]\Big] \\
&&\dis\le  C\Big(\dbE^{\dbP_0}\big[ (M_T^{\a^\e})^{-2}\big]\Big)^{1\over 8} \Big(\dbE^{\dbP_0}\big[M_T^{\a^\e}\int_0^T \D f_i(s, B,  \a^\e_s, \hat Z^{\e, i}_s)ds\big]\Big)^{1\over 4}   \Big(\dbE^{\dbP_0}\big[ \int_0^T [1+|\hat Z^{\e, i}_s|2]\big]\Big)^{5\over 8}\\
&&\dis\le C\e^{1\over 4}. 
\eeaa
This, together with \reff{xiest}, implies that 
\beaa
\dbE^{\dbP_0}\Big[ |\xi_i(B) - Y^{y, \a^{\e}, \hat Z^{\e}, i}_T|^2  + \int_0^T  [\D f_i(s, B,  \a^\e_s, \hat Z^{\e, i}_s)]^{3\over 2} ds\Big] \le C\e^{1\over 4}.
\eeaa
Then, by \reff{W} we have $W(0,0, y) \le CN\e^{1\over 4}$. Since $\e$ is arbitrary, we get $W(0,0,y) = 0$, that is,  $y\in \dbN(0,0)$.
\qed

Note that \reff{W} is a standard path dependent control problem. Following Zhang \cite[Section 11.3.3]{Zhang} we have the following result whose proof is omitted.
 \begin{prop}
 \label{prop-W}
 Let Assumptions \ref{assum-bound} and \ref{assum-reg}  hold. Then  $W\in C([0, T\times \O \times \dbR^N)$ is a viscosity solution of the following path dependent PDE:
 \bea
 \label{PPDE}
 \left.\ba{c}
 \dis \pa_t W + \inf_{a\in A, z\in \dbR^{N\times d}}\Big[{1\over 2} \tr(\pa^2_{\o\o} W) + {1\over 2} \tr(z^\top  \pa^2_{y y} W z)+ \tr( z^\top  \pa_{y \o} W)\\
  \dis \qquad + \sum_{i=1}^N \big[ [\D f_i(t, \o, a, z_i) ]^{3\over 2}-\underline f_i(t,\o, a^{-i}, z_i) \pa_{y_i} W \big]\Big]  =0; \\
\dis  W(T,\o, y) = |\xi(\o) - y|^2.
 \ea\right.
 \eea 
 \end{prop}

\begin{rem}
\label{rem-PPDE}
{\rm \begin{enumerate}
\item The path derivatives $\pa_\o W, \pa^2_{\o\o} W$ are introduced by Dupire \cite{Dupire}, and we refer to Zhang \cite[Section 9.4]{Zhang} for more details. Note that this path dependent PDE is always degenerate and the control is unbounded,  so the uniqueness of viscosity solution is not completely covered by Ekren, Touzi, \& Zhang \cite{ETZ1, ETZ2} and Ren, Touzi, \& Zhang \cite{RTZ}. This problem is in general challenging and is left for future research.

\item In the state dependent case as in Remark \ref{rem-cont-state}, $W = W(t,x,y)$ also becomes state dependent and the path dependent PDE \reff{PPDE} reduces to a standard HJB equation:
 \beaa
 \label{PDE}
 \left.\ba{c}
 \dis \pa_t W + \inf_{a\in A, z\in \dbR^{N\times d}}\Big[{1\over 2} \tr(\pa^2_{xx} W) + {1\over 2} \tr(z^\top  \pa^2_{y y} W z)+ \tr( z^\top  \pa_{y x} W)\\
  \dis + \sum_{i=1}^N \big[ [\D f_i(t, x, a, z_i) ]^{3\over 2}-\underline f_i(t,x, a^{-i}, z_i) \pa_{y_i} W \big]\Big]  =0; \\
\dis  W(T,x, y) = |g(x) - y|^2.
 \ea\right.
 \eeaa 
 This PDE is also degenerate  and with unbounded controls though.
 
 {\color{black}
 \item In light of Theorem \ref{thm-duality},   PPDE \reff{PPDE}, especially  PDE \reff{PDE} in the state dependent case,  is quite useful for numerical computation of the set value $\dbV(t,\o)$. }
\qed
\end{enumerate}}
\end{rem}

\begin{rem}
\label{rem-viability}
{\rm Roughly speaking (modulus the existence of optimal controls in \reff{W}), $y$ is in the nodal set $N(t, \o)$ if and only if there exists $\a, Z$ such that $Y^{t,\o, y, \a, Z}$ in \reff{X} hits the target $\xi^{t,\o}(B)$ at $T-t$.  This is in the spirit of Cardaliaguet,  Quincampoix, \&  Saint-Pierre \cite{CQS}. However, we note that \cite{CQS} uses strategy versus controls, while we use closed-loop controls for all players.  
\qed}
\end{rem}

{\color{black}  
\begin{rem}
\label{rem-BSDE}
{\rm In this remark we make further connection between the game and BSDEs.

\begin{enumerate}

\item{} In the literature, one may indeed use \reff{BSDENE} to find equilibriums, especially in the state dependent setting \reff{cont-state}, see e.g. Hamadene, Lepeltier, \& Peng \cite{HLP},  Hamadene \& Mu \cite{HM1, HM2}, and Espinosa \& Touzi \cite{ET}. To be precise, assume there exist a measurable function $\f: [0, T]\times \dbR^d\times (\dbR^d)^N \to A$ such that, for $i=1,\cds,N$,
\bea
\label{NEphi}
\underline f_i(t, x, \f^{-i}(t,x,z), z^i) =   f_i(t, x,  \f(t,x,z), z^i),
\eea
and the following BSDEs have a strong solution (setting $(t,x)=(0,0)$ for simplicity):
\bea
\label{BSDEphi}
 \underline Y^{i}_s = g_i(B_{T}) +  \int_s^{T} \underline f_i(r, B_r, \f^{-i}(r, B_r, \underline Z_r), \underline Z^{i}_r) dr - \int_s^{T}  \underline Z^{i}_r dB_r,
 \eea
then $\a^*_t := \f(t, B_t, \underline Z_t)$ is a Nash equilibrium at $(0,0)$. However, we should note that the function $\f$, assuming its existence,  may not be continuous  and thus the wellposedness of \reff{BSDEphi} may not be easy. Even worse, in order to obtain the whole set $\dbV_0(0,0)$, as we noted before we need to consider path dependent $\f: [0, T]\times \O\times (\dbR^d)^N \to A$,  which will make the wellposedness of \reff{BSDEphi} even harder. Nevertheless, by \reff{BSDENE} it is true that the set $\dbV_0$ can be constructed by first finding all path dependent functions $\f$ satisfying  \reff{NEphi} and then finding all strong solutions of the multidimensional BSDE \reff{BSDEphi}, where both \reff{NEphi} and \reff{BSDEphi} should be extended to the path dependent setting.

\item{} One may replace the linear BSDE \reff{BSDE} with nonlinear BSDEs: 
\beaa
Y^{t,\o, \a,i}_s = \xi^{t,\o}_i(B) + \int_s^{T-t} f_i^{t,\o}(r, B, \a_r, Y^{ t,\o,\a,i}_r, Z^{ t,\o,\a,i}_r) dr -\!\! \int_s^{T-t} Z^{t,\o, \a, i}_r dB_r,
\eeaa 
where $f_i: [0, T]\times \O \times A\times \dbR\times \dbR^d \to \dbR$ is nonlinear in $(y,z)$. Still define $J(t,\o,\a) := Y^{t,\o,\a}_0$, then one can show without significant difficulties that all the results in this section hold true after obvious modifications. 
\qed
\end{enumerate}
}
\end{rem} 
}

\end{document}